\documentclass[12pt,leqno]{article}
\usepackage{amssymb}
\usepackage[mathscr]{eucal}
\usepackage{amsmath,amssymb,latexsym,theorem,bbm} %
\usepackage{color,url}
\usepackage{appendix}

\newcommand{\NN}{\mathbb{N}}

\newcommand{\RR}{\mathbb{R}}

\newcommand{\ZZ}{\mathbb{Z}}

\newcommand{\be}{{\boldsymbol{e}}}

\newcommand{\tL}{\widetilde{L}}

\newcommand{\tm}{{\widetilde{m}}}
\newcommand{\bm}{{\boldsymbol{m}}}
\newcommand{\bM}{{\boldsymbol{M}}}

\newcommand{\bV}{{\boldsymbol{V}}}

\newcommand{\bX}{{\boldsymbol{X}}}

\newcommand{\bY}{{\boldsymbol{Y}}}

\newcommand{\bxi}{{\boldsymbol{\xi}}}

\newcommand{\bvare}{{\boldsymbol{\vare}}}
\newcommand{\bpi}{{\boldsymbol{\pi}}}

\newcommand{\bzero}{{\boldsymbol{0}}}

\newcommand{\cA}{{\mathcal A}}

\newcommand{\cD}{{\mathcal D}}

\newcommand{\cF}{{\mathcal F}}

\newcommand{\dd}{\mathrm{d}}
\newcommand{\ee}{\mathrm{e}}

\newcommand{\SUB}{{\mathrm{sub}}}
\newcommand{\CRIT}{{\mathrm{crit}}}
\newcommand{\SUP}{{\mathrm{sup}}}

\newcommand{\INARp}{\textup{INAR($p$)}}

\newcommand{\EE}{\operatorname{\mathbb{E}}}
\newcommand{\PP}{\operatorname{\mathbb{P}}}
\newcommand{\OO}{\operatorname{O}}
\newcommand{\oo}{\operatorname{o}}
\newcommand{\var}{\operatorname{Var}}

\newcommand{\tM}{\widetilde{M}}

\newcommand{\tS}{\widetilde{S}}

\newcommand{\ty}{\widetilde{y}}

\newcommand{\vare}{\varepsilon}
\newcommand{\tzeta}{\widetilde{\zeta}}

\renewcommand{\mid}{\,|\,}

\renewcommand{\leq}{\leqslant}
\renewcommand{\geq}{\geqslant}

\newcommand{\distr}{\stackrel{\cD}{\longrightarrow}}

\newcommand{\distre}{\stackrel{\cD}{=}}

\newcommand{\as}{\stackrel{{\mathrm{a.s.}}}{\longrightarrow}}

\newcommand{\bbone}{\mathbbm{1}}

\newcommand{\proofend}{\hfill\mbox{$\Box$}}

\numberwithin{equation}{section}

\theoremstyle{change} \theorembodyfont{\em}
\newtheorem{Thm}{Theorem.}[section]
\newtheorem{Lem}[Thm]{Lemma.}
\newtheorem{Pro}[Thm]{Proposition.}

\newtheorem{Def}[Thm]{Definition.}

\theorembodyfont{\rm}
\newtheorem{Rem}[Thm]{Remark.}

\begin{document}

\begin{center}
 {\bfseries\large On tail behaviour of stationary second-order Galton--Watson processes with
                   immigration} \\[6mm]
 {\sc M\'aty\'as $\text{Barczy}^{*,\diamond}$,
            \ Zsuzsanna $\text{B\H{o}sze}^{**}$,
            \ Gyula $\text{Pap}^{***}$}
\end{center}

\vskip0.2cm

\noindent
 * MTA-SZTE Analysis and Stochastics Research Group,
   Bolyai Institute, University of Szeged,
   Aradi v\'ertan\'uk tere 1, H--6720 Szeged, Hungary.

\noindent
 ** Institute for Mathematical Stochastics,
    Georg-August-Universit\"at
    G\"ottingen, Goldschmidtstr. 7, 37077 G\"ottingen, Germany.

\noindent
 *** Bolyai Institute, University of Szeged,
    Aradi v\'ertan\'uk tere 1, H-6720 Szeged, Hungary.

\noindent e-mail: barczy@math.u-szeged.hu (M. Barczy),
                  zsuzsanna.boesze@uni-goettingen.de (Zs. B\H{o}sze).

\noindent $\diamond$ Corresponding author.



\renewcommand{\thefootnote}{}
\footnote{\textit{2010 Mathematics Subject Classifications\/}:
 60J80, 60G70.}
\footnote{\textit{Key words and phrases\/}:
 second-order Galton--Watson process with immigration, regularly varying distribution, tail
 behavior.}
\vspace*{0.2cm}
\footnote{
Supported by the Hungarian Croatian Intergovernmental S \& T Cooperation Programme for 2017-2018
 under Grant No.\ 16-1-2016-0027.
M\'aty\'as Barczy is supported by the J\'anos Bolyai Research Scholarship of the Hungarian
 Academy of Sciences.}

\vspace*{-10mm}

\begin{abstract}
A second-order Galton--Watson process with immigration can be represented
 as a coordinate process of a 2-type Galton--Watson process with  immigration.
Sufficient conditions are derived on the offspring and immigration
 distributions of a second-order Galton--Watson process with immigration under which
 the corresponding 2-type Galton--Watson process with immigration has a unique stationary
 distribution such that its common marginals are regularly varying.
In the course of the proof sufficient conditions are given under which the distribution of a
 second-order Galton--Watson process (without immigration) at any fixed time is regularly varying
 provided that the initial sizes of the population are independent and regularly varying.
\end{abstract}

\section{Introduction}
\label{section_intro}

Branching processes have been frequently used in biology, e.g., for modeling the spread of an
 infectious disease,
 for gene amplification and deamplification or for modeling telomere shortening, see, e.g., Kimmel and
 Axelrod \cite{KimAxe}.
Higher-order Galton--Watson processes with immigration having finite  second moment
 (also called Generalized Integer-valued AutoRegressive (GINAR)  processes) have been
 introduced by Latour \cite[equation (1.1)]{Lat}.
P\'enisson and Jacob \cite{PenJac} used higher-order Galton--Watson processes (without immigration) for studying the decay phase of an epidemic, and, as an application, they investigated the Bovine Spongiform Encephalopathy epidemic in Great Britain after the 1988 feed ban law.
As a continuation, P\'enisson \cite{Pen} introduced estimators of the so-called infection parameter in the growth
 and decay phases of an epidemic.
Recently, Kashikar and Deshmukh \cite{KasDes, KasDes2} and  Kashikar \cite{Kas} used second order Galton--Watson processes (without immigration) for modeling the swine flu data for Pune, India and La-Gloria, Mexico.
Kashikar and Deshmukh \cite{KasDes} also studied their basic probabilistic properties such as a formula for their probability generator function, probability of extinction, long run behavior and conditional least squares
 estimation of the offspring means.
Higher-order Galton--Watson processes with immigration are special multi-type Galton--Watson processes with immigration,
 and to give an example for an application of such processes for modeling epidemics, for example, we can mention
 D\'enes et al.\ \cite{DenKevNisRos}, where a 17-type Galton--Watson process with immigration has been applied to
 describe the risk of a major epidemic in connection with the 2012 UEFA European Football Championship
 took place in Ukraine and Poland between 8 June and 1 July 2012.

Let \ $\ZZ_+$, \ $\NN$, \ $\RR$, \ $\RR_+$, \ $\RR_{++}$, \ and \ $\RR_{--}$ \ denote the set of
 non-negative integers, positive integers, real numbers, non-negative real numbers, positive real
 numbers and negative real numbers, respectively.
For functions \ $f : \RR_{++} \to \RR_{++}$ \ and \ $g : \RR_{++} \to \RR_{++}$, \  by the
 notation \ $f(x) \sim g(x)$, \ $f(x) = \oo(g(x))$ \ and \ $f(x) = \OO(g(x))$ \ as
 \ $x \to \infty$, \ we mean that \ $\lim_{x\to\infty} \frac{f(x)}{g(x)} = 1$,
 \ $\lim_{x\to\infty} \frac{f(x)}{g(x)} = 0$ \ and
 \ $\limsup_{x\to\infty} \frac{f(x)}{g(x)} < \infty$, \ respectively.
The natural basis of \ $\RR^d$ \ will be denoted by \ $\{\be_1, \ldots, \be_d\}$.
\ For \ $x \in \RR$, \ the integer part of \ $x$ \ is denoted by \ $\lfloor x\rfloor$.
\ Every random variable will be defined on a probability space \ $(\Omega, \cA, \PP)$.
\ Equality in distributions of random variables or stochastic processes is denoted by \ $\distre$.

First, we recall the Galton--Watson process with immigration, which assumes that an individual can
 reproduce only once during its lifetime at age 1, and then it dies immediately.
The initial population size at time \ $0$ \ will be denoted by \ $X_0$.
\ For each \ $n \in \NN$, \ the population consists of the offsprings born at time \ $n$ \ and the
 immigrants arriving at time \ $n$.
\ For each \ $n, i \in \NN$, \ the number of offsprings produced at time \ $n$ \ by the
 \ $i^\mathrm{th}$ \ individual of the \ $(n-1)^\mathrm{th}$ \ generation will be denoted by
 \ $\xi_{n,i}$.
\ The number of immigrants in the \ $n^\mathrm{th}$ \ generation will be denoted by \ $\vare_n$.
\ Then, for the population size \ $X_n$ \ of the \ $n^\mathrm{th}$ \ generation, we have
 \begin{equation}\label{GWI}
  X_n = \sum_{i=1}^{X_{n-1}} \xi_{n,i} + \vare_n , \qquad n \in \NN ,
 \end{equation}
  where \ $\sum_{i=1}^0:=0$.
\ Here \ $\bigl\{X_0, \, \xi_{n,i}, \, \vare_n : n, i \in \NN\bigr\}$ \ are supposed to be
 independent non-negative integer-valued random variables, and
 \ $\{\xi_{n,i} : n, i \in \NN\}$ \ and \ $\{\vare_n : n \in \NN\}$ \ are supposed to consist of
 identically distributed random variables, respectively.
If \ $\vare_n = 0$, \ $n \in \NN$, \ then we say that \ $(X_n)_{n\in\ZZ_+}$ \ is a Galton--Watson
 process (without immigration).

Next, we introduce the second-order Galton--Watson branching model with immigration.
In this model we suppose that an individual reproduces at age \ $1$ \ and also at age \ $2$, \ and
 then it dies immediately.
For each \ $n \in \NN$, \ the population consists again of the offsprings born at time \ $n$ \ and
 the immigrants arriving at time \ $n$.
\ For each \ $n, i, j \in \NN$, \ the number of offsprings produced at time \ $n$ \ by the
 \ $i^\mathrm{th}$ \ individual of the \ $(n-1)^\mathrm{th}$ \ generation and by the
 \ $j^\mathrm{th}$ \ individual of the \ $(n-2)^\mathrm{nd}$ \ generation will be denoted by
 \ $\xi_{n,i}$ \ and \ $\eta_{n,j}$, \ respectively, and \ $\vare_n$ \ denotes the number of
 immigrants in the \ $n^\mathrm{th}$ \ generation.
Then, for the population size \ $X_n$ \ of the \ $n^\mathrm{th}$ \ generation, we have
 \begin{equation}\label{2oGWI}
  X_n = \sum_{i=1}^{X_{n-1}} \xi_{n,i} + \sum_{j=1}^{X_{n-2}} \eta_{n,j} + \vare_n ,
  \qquad n \in \NN ,
 \end{equation}
 where \ $X_{-1}$ \ and \ $X_0$ \ are non-negative integer-valued random variables (the initial
 population sizes).
Here \ $\bigl\{ X_{-1}, X_0, \, \xi_{n,i}, \, \eta_{n,j}, \, \vare_n : n, i, j \in \NN \bigr\}$
 \ are supposed to be non-negative integer-valued random variables such that
 \ $\bigl\{ (X_{-1}, X_0), \, \xi_{n,i}, \, \eta_{n,j}, \, \vare_n : n, i, j \in \NN \bigr\}$ \ are independent,
 and \ $\{\xi_{n,i} : n, i \in \NN\}$, \ $\{\eta_{n,j} : n, j \in \NN\}$ \ and
 \ $\{ \vare_n : n \in \NN\}$ \ are supposed to consist of identically distributed random
 variables, respectively.
Note that the number of individuals alive at time \ $n \in \ZZ_+$ \ is \ $X_n + X_{n-1}$, \ which
 can be larger than the population size \ $X_n$ \ of the \ $n^\mathrm{th}$ \ generation, since the
 individuals of the population at time \ $n-1$ \ are still alive at time \ $n$, \ because they can
 reproduce also at age \ $2$.
\ The stochastic process \ $(X_n)_{n\geq -1}$ \ given by \eqref{2oGWI} is called a second-order
 Galton--Watson process with immigration or a Generalized Integer-valued AutoRegressive process
 of order 2 (GINAR(2) process), see, e.g., Latour \cite{Lat}.
Especially, if \ $\xi_{1,1}$ \ and \ $\eta_{1,1}$ \ are Bernoulli distributed random variables,
 then \ $(X_n)_{n\geq -1}$ \ is also called an Integer-valued AutoRegressive process of order 2
 (INAR(2) process), see, e.g., Du and Li \cite{DuLi}.
If \ $\vare_1 = 0$, \ then we say that \ $(X_n)_{n\geq-1}$ \ is a second-order Galton--Watson
 process without immigration, introduced and studied by Kashikar and Deshmukh \cite{KasDes} as well.

The process given in \eqref{2oGWI} with the special choice \ $\eta_{1,1} = 0$ \ gives back the
 process given in \eqref{GWI}, which will be called a first-order Galton--Watson process with
 immigration to make a distinction.

For notational convenience, let \ $\xi$, \ $\eta$ \ and \ $\vare$ \ be random variables such that
 \ $\xi \distre \xi_{1,1}$, \ $\eta \distre \eta_{1,1}$ \ and \ $\vare \distre \vare_1$,
 \ and put \ $m_\xi := \EE(\xi) \in [0, \infty]$, \ $m_\eta := \EE(\eta) \in [0, \infty]$
 \ and \ $m_\vare := \EE(\vare) \in [0, \infty]$.

If \ $(X_n)_{n\in\ZZ_+}$ \ is a (first-order) Galton--Watson process with immigration such that
 \ $m_\xi \in (0, 1)$, \ $\PP(\vare=0)<1$, \  and \ $\sum_{j=1}^\infty \PP(\vare = j) \log(j) < \infty$, \ then the
 Markov process \ $(X_n)_{n\in\ZZ_+}$ \ admits a unique stationary distribution \ $\mu$, \ see,
 e.g., Quine \cite{Qui}.
If \ $\vare$ \ is regularly varying with index \ $\alpha \in \RR_{++}$, \ i.e.,
 \ $\PP(\vare > x) \in \RR_{++}$ \ for all \ $x \in \RR_{++}$, \ and
 \[
 \lim_{x\to\infty}
  \frac{\PP(\vare>qx)}{\PP(\vare>x)} = q^{-\alpha} \qquad \text{for all \ $q \in \RR_{++}$,}
 \]
 then, by Lemma \ref{help_log_exp}, \ $\sum_{j=1}^\infty \PP(\vare = j) \log(j) < \infty$.
\ The content of Theorem 2.1.1 in Basrak et al.\ \cite{BasKulPal} is the following statement.

\begin{Thm}\label{GWI_stat}
Let \ $(X_n)_{n\in\ZZ_+}$ \ be a (first-order) Galton--Watson process with immigration such that
 \ $m_\xi \in (0, 1)$ \ and \ $\vare$ \ is regularly varying with index \ $\alpha \in (0, 2)$.
\ In case of \ $\alpha \in [1, 2)$, \ assume additionally that \ $\EE(\xi^2) < \infty$.
\ Then the tail of the unique stationary distribution \ $\mu$ \ of
 \ $(X_n)_{n\in\ZZ_+}$ \ satisfies
 \[
   \mu((x,\infty)) \sim \sum_{i=0}^\infty m_\xi^{i\alpha}\,\PP(\vare > x)
                         = \frac{1}{1-m_\xi^\alpha} \PP(\vare > x) \qquad
   \text{as \ $x \to \infty$,}
 \]
 and hence \ $\mu$ \ is also regularly varying with index \ $\alpha$.
\end{Thm}

Note that in case of \ $\alpha = 1$ \ and \ $m_\vare = \infty$ \ Basrak et al.\
 \cite[Theorem 2.1.1]{BasKulPal} assume additionally that \ $\vare$ \ is consistently varying (or
 in other words intermediate varying), but, eventually, it follows from the fact that \ $\vare$ \ is regularly
 varying.
Basrak et al.\ \cite[Remark 2.2.2]{BasKulPal} derived the result of Theorem \ref{GWI_stat} also for
 \ $\alpha \in [2, 3)$ \ under the additional assumption \ $\EE(\xi^3) < \infty$ \ (not mentioned
 in the paper), and they remark that the same applies to all \ $\alpha \in [3, \infty)$
 \ (possibly under an additional moment assumption
 \ $\EE(\xi^{\lfloor \alpha\rfloor + 1})<\infty$).

In Barczy et al.\ \cite{BarBosPap} we study regularly varying non-stationary (first-order) Galton--Watson processes with immigration.

As the main result of the paper, in Theorem \ref{Thm3}, in the same spirit as in Theorem \ref{GWI_stat},
 we present sufficient conditions on the offspring and immigration distributions of a second-order
 Galton--Watson process with immigration under which its associated 2-type Galton--Watson process
 with immigration has a unique stationary distribution such that its common marginals are regularly varying.
According to our knowledge, such a result has not been established so far,
 e.g., we could not find any reference which would address regularly varying GINAR(2) processes.
 Our result and the applied technique might be extended to a \ $p$-th order Galton--Watson branching process with immigration,
 however such an extension is not immediate, for example, it is not clear what would replace the constant
 \ $\sum_{i=0}^\infty m_i^\alpha$ \ in Theorem \ref{Thm3}.
\ More generally, one can pose an open problem, namely, under what conditions on the offspring and immigration distributions
 of a general \ $p$-type Galton--Watson branching process with immigration, its unique ($p$-dimensional) stationary distribution
  is jointly regularly varying.
We also note that there is a vast literature on tail behavior of regularly varying time series (see, e.g., Hult and Samorodnitsky \cite{HulSam}),
 however, the available results do not seem to be applicable for describing the tail behavior of the stationary distribution for regularly varying
 branching processes.
The link between GINAR and autoregressive processes is that their autocovariance functions are identical under finite second moment assumptions,
 but we can not see that it would imply anything for the tail behavior of a GINAR process knowing the tail behaviour of a corresponding
 autoregressive process.
Further, in our situation the second moment is infinite, so the autocovariance function is not defined.

Very recently, B\H osze and Pap \cite{BosPap} have studied regularly varying non-stationary second-order
 Galton--Watson processes with immigration.
They have found some sufficient conditions on the initial, the offspring
 and the immigration distributions of a non-stationary second-order Galton–-Watson process with immigration
 under which the distribution of the process in question is regularly varying at any fixed time.
The results in B\H osze and Pap \cite{BosPap} can be considered as extensions of the results
 in Barczy et al.\ \cite{BarBosPap} on not necessarily stationary (first-order) Galton--Watson processes with immigration.
Concerning the results in B\H osze and Pap \cite{BosPap} and in the present paper, there is no overlap,
 for more details see Remark \ref{Rem_no_overlap}.

The paper is organized as follows.
In Section \ref{section_2GWI}, first, for a second-order  Galton--Watson process with immigration,
 we give a representation of the unique stationary distribution and its marginals, respectively,
 then our main result, Theorem \ref{Thm3}, is formulated.
The rest of Section \ref{section_2GWI} is devoted to the proof of Theorem \ref{Thm3}.
In the course of the proof, we formulate an auxiliary result about the tail behaviour of a second-order Galton--Watson process
 (without immigration) with a regularly varying initial distribution at time \ $0$ \ and with value \ $0$ \ at time \ $-1$,
 \ see Proposition \ref{Cor_2GW_X_0_X_-1}.
We close the paper with seven appendices which are used throughout the proofs.
In Appendix \ref{section_prel}, we recall a representation of a second-order Galton--Watson process without or with
 immigration as a (special) 2-type Galton--Watson process without or with immigration, respectively.
In Appendix \ref{App0}, we derive an explicit formula for the expectation of a second-order Galton--Watson process with immigration at
 time \ $n$ \ and describe its asymptotic behavior as \ $n \to \infty$.
Appendix \ref{App_moments} is about the existence and estimation of higher order moments of a second-order Galton--Watson process (without immigration).
In Appendix \ref{App2multitype}, we recall a representation of the unique stationary distribution for a 2-type Galton--Watson
 process with immigration.
In Appendix \ref{App1}, we collect several results on regularly varying functions and distributions, to name a few of them:
 convolution property, Karamata's theorem and Potter's bounds.
Appendix \ref{App3} is devoted to recall and (re)prove a result on large deviations for sums of non-negative independent
 and identically distributed regularly varying random variables due to Tang and Yan \cite[part (ii) of Theorem 1]{TangYan}.
Finally, in Appendix \ref{section_2GW}, we present a variant of Proposition \ref{Cor_2GW_X_0_X_-1},
 where the initial values \ $X_{-1}$ \ and \ $X_0$ \ are independent and regularly varying
 together with a second type of proof, see Proposition \ref{2GW_X_0_X_-1}.

\section{Tail behavior of the marginals of the stationary distribution of second-order Galton--Watson
          processes with immigration}
\label{section_2GWI}

Let \ $(X_n)_{n\geq-1}$ \ be a second order Galton--Watson process with immigration given in
 \eqref{2oGWI}, and let us consider the Markov chain \ $(\bY_k)_{k\in\ZZ_+}$ \ given by
 \begin{align*}
  \bY_n := \begin{bmatrix}
            Y_{n,1} \\
            Y_{n,2} \\
           \end{bmatrix}
        := \begin{bmatrix}
            X_n \\
            X_{n-1} \\
           \end{bmatrix}
        =\sum_{i=1}^{Y_{n-1,1}}
           \begin{bmatrix} \xi_{n,i} \\ 1 \end{bmatrix}
          + \sum_{j=1}^{Y_{n-1,2}}
             \begin{bmatrix} \eta_{n,j} \\ 0 \end{bmatrix}
          + \begin{bmatrix} \vare_n \\ 0 \end{bmatrix} ,
  \qquad n \in \NN ,
 \end{align*}
 which is a (special) 2-type Galton--Watson process with immigration, and
 \ $(\be_1^\top \bY_k)_{k\in\ZZ_+} = (X_k)_{k\in\ZZ_+}$,
 \ $(\be_2^\top \bY_{k+1})_{k\geq-1} = (X_k)_{k\geq-1}$ \
 (for more details, see Appendix \ref{section_prel}).
\ If \ $m_\xi \in \RR_{++}$, \ $m_\eta \in \RR_{++}$, \ $m_\xi + m_\eta < 1$,
 \ $\PP(\vare = 0) < 1$ \ and \ $\EE(\bbone_{\{\vare\ne0\}} \log(\vare)) < \infty$, \ then there
 exists a unique stationary distribution \ $\bpi$ \ for \ $(\bY_n)_{n\in\ZZ_+}$, \ see Appendix
 \ref{App2multitype}, since then \ $\bM_{\xi, \eta}$ \ is primitive due to the fact that
 \[
   \bM_{\xi, \eta}^2
   = \begin{bmatrix}
      m_\xi & m_\eta \\
      1 & 0
     \end{bmatrix}^2
   = \begin{bmatrix}
      m_\xi^2+m_\eta & m_\xi m_\eta \\
      m_\xi & m_\eta
    \end{bmatrix}
   \in \RR_{++}^2 .
 \]
Moreover, the stationary distribution \ $\bpi$ \ of \ $(\bY_n)_{n\in\ZZ_+}$ \ has a representation
 \[
   \bpi \distre \sum_{i=0}^\infty \bV_i^{(i)}(\bvare_i) ,
 \]
 where \ $(\bV_k^{(i)}(\bvare_i))_{k\in\ZZ_+}$, \ $i \in \ZZ_+$, \ are independent copies of a
 (special) 2-type Galton--Watson process \ $(\bV_k(\bvare))_{k\in\ZZ_+}$ \ (without immigration)
 with initial vector \ $\bV_0(\bvare) = \bvare$ \ and with the same offspring distributions as
 \ $(\bY_k)_{k\in\ZZ_+}$, \ and the series \ $\sum_{i=0}^\infty \bV_i^{(i)}(\bvare)$ \ converges
 with probability 1, see Appendix \ref{App2multitype}.
Using the considerations for the backward representation in Appendix \ref{section_prel}, we have
 \ $(\be_1^\top \bV_k(\bvare))_{k\in\ZZ_+} = (V_k(\vare))_{k\in\ZZ_+}$ \ and
 \ $(\be_2^\top \bV_{k+1}(\bvare))_{k\geq-1} = (V_k(\vare))_{k\geq-1}$, \ where
 \ $(V_k(\vare))_{k\geq-1}$ \ is a second-order Galton--Watson process (without immigration) with
 initial values \ $V_0(\vare) = \vare$ \ and \ $V_{-1}(\vare) = 0$, \ and with the same offspring
 distributions as \ $(X_k)_{k\geq-1}$.
\ Consequently, the marginals of the stationary distribution \ $\bpi$ \ are the same distributions
 \ $\pi$, \ and it admits the representation
 \[
   \pi \distre \sum_{i=0}^\infty V_i^{(i)}(\vare_i) ,
 \]
 where \ $(V_k^{(i)}(\vare_i))_{k\in\ZZ_+}$, \ $i \in \ZZ_+$, \ are independent copies of
 \ $(V_k(\vare))_{k\geq-1}$.
\ This follows also from the fact that the stationary distribution \ $\bpi$ \ is the limit in
 distribution of \ $ \bY_n$ \ as \ $n \to \infty$ \ and
 \[
   \bY_n
   = \begin{bmatrix}
      X_n \\
      X_{n-1}
     \end{bmatrix} ,
   \qquad n \in \ZZ_+ ,
 \]
 thus the coordinates of \ $ \bY_n$ \ converge in distribution to the same distribution \ $\pi$
 \ as \ $n \to \infty$.

Note that \ $(X_n)_{n\geq-1}$ \ is only a second-order Markov chain, but not a Markov chain.
Moreover, \ $(X_n)_{n\geq-1}$ \ is strictly stationary if and only if the distribution of the
 initial population sizes \ $(X_0, X_{-1})^\top$ \ coincides with the stationary distribution
 \ $\bpi$ \ of the Markov chain \ $(\bY_k)_{k\in\ZZ_+}$.
\ Indeed, if \ $(X_0, X_{-1})^\top \distre \bpi$, \ then \ $\bY_0 \distre \bpi$, \ thus
 \ $(\bY_k)_{k\in\ZZ_+}$ \ is strictly stationary, and hence for each \ $n, m \in \ZZ_0$,
 \ $(\bY_0, \ldots, \bY_n) \distre (\bY_m, \ldots, \bY_{n+m})$, \ yielding
 \[
   (X_0, X_{-1}, X_1, X_0, \ldots, X_n, X_{n-1})
   \distre (X_m, X_{m-1}, X_{m+1}, X_m, \ldots, X_{n+m}, X_{n+m-1}) .
 \]
Especially, \ $(X_{-1}, X_0, X_1, \ldots, X_n) \distre (X_{m-1}, X_m, X_{m+1}, \ldots, X_{n+m})$,
 \ hence \ $(X_n)_{n\geq-1}$ \ is strictly stationary.
Since \ $(X_m, X_{m-1}, X_{m+1}, X_m, \ldots, X_{n+m}, X_{n+m-1})$ \ is a continuous function of
 \ $(X_{m-1}, X_m, X_{m+1}, \ldots, X_{n+m})$, \ these considerations work backwards as well.
Consequently, \ $\bpi$ \ is the unique stationary distribution of the second-order Markov chain
 \ $(X_n)_{n\geq-1}$.

\begin{Thm}\label{Thm3}
Let \ $(X_n)_{n\geq-1}$ \ be a second-order Galton--Watson process with immigration such that
 \ $m_\xi \in \RR_{++}$, \ $m_\eta \in \RR_{++}$, \ $m_\xi + m_\eta < 1$ \ and \ $\vare$ \ is
 regularly varying with index \ $\alpha \in (0, 2)$.
\ In case of \ $\alpha \in [1, 2)$, \ assume additionally that \ $\EE(\xi^2) < \infty$ \ and
 \ $\EE(\eta^2) < \infty$.
\ Then the tail of the marginals \ $\pi$ \ of the unique stationary distribution \ $\bpi$ \ of
 \ $(X_n)_{n\geq-1}$ \ satisfies
 \[
   \pi((x,\infty)) \sim \sum_{i=0}^\infty m_i^\alpha \,\PP(\vare > x) \qquad
   \text{as \ $x \to \infty$,}
 \]
 where \ $m_0 := 1$ \ and
 \begin{align}\label{m_n}
  m_k := \frac{\lambda_+^{k+1}-\lambda_-^{k+1}}{\lambda_+-\lambda_-} , \qquad
  \lambda_+ := \frac{m_\xi+\sqrt{m_\xi^2+4m_\eta}}{2} , \qquad
  \lambda_- := \frac{m_\xi-\sqrt{m_\xi^2+4m_\eta}}{2}
 \end{align}
 for \ $k \in \NN$.
\ Consequently, \ $\pi$ \ is also regularly varying with index \ $\alpha$.
\end{Thm}

Note that \ $\lambda_+$ \ and \ $\lambda_-$ \ are the eigenvalues of the offspring mean
 matrix \ $\bM_{\xi,\eta}$ \ given in \eqref{bA} related to the recursive formula \eqref{recEX_n}
 for the expectations \ $\EE(X_n)$, \ $n \in \NN$.
\ For each \ $k \in \ZZ_+$, \ the assumptions \ $m_\xi \in \RR_{++}$ \ and \ $m_\eta \in \RR_{++}$
 \ imply \ $m_k \in \RR_{++}$.
\ Further, by \eqref{EXn}, for all \ $k \in \ZZ_+$, \ we have \ $m_k = \EE(V_{k,0})$, \ where
 \ $(V_{n,0})_{n\geq-1}$ \ is a second-order Galton--Watson process (without immigration) with
 initial values \ $V_{0,0} = 1$ \ and \ $V_{-1,0} = 0$, \ and with the same offspring distributions
 as \ $(X_n)_{n\geq-1}$.
Consequently, the series \ $\sum_{i=0}^\infty m_i^{\alpha}$ \ appearing in Theorem \ref{Thm3} is convergent,
 since for each \ $i \in \NN$, \ we have \ $m_i = \EE(V_{i,0}) \leq \lambda_+^i < 1$ \ by \eqref{1moment_est} and
 the assumption \ $m_\xi + m_\eta < 1$.

We point out that in Theorem \ref{Thm3} only the regular variation of the marginals \ $\pi$ \ of
 \ $\bpi$ \ is proved, the question of the joint regular variation of \ $\bpi$ \ remains open.

\begin{Rem}\label{Rem_no_overlap}
Note that there is no overlap between the results in the recent paper of B\H osze and Pap \cite{BosPap}
 on non-stationary second-order Galton–-Watson processes with immigration
 and in the present paper.
In \cite{BosPap} the authors always suppose that the initial values \ $X_0$ \ and \ $X_{-1}$ \ of
 a second-order Galton–-Watson process with immigration \ $(X_n)_{n\geq -1}$ \ are independent,
 so in the results of \cite{BosPap} the distribution of \ $(X_0,X_{-1})$ \ can not be chosen as
 the unique stationary distribution \ $\bpi$, \ since the marginals of \ $\bpi$ \ are not independent in general.
\proofend
\end{Rem}

For the proof of Theorem \ref{Thm3}, we need an auxiliary result on the tail behaviour of
 second-order Galton--Watson processes (without immigration) having regularly varying initial
 distributions.

\begin{Pro}\label{Cor_2GW_X_0_X_-1}
Let \ $(X_n)_{n\geq-1}$ \ be a second-order Galton--Watson process (without immigration) such
 that \ $X_0$ \ is regularly varying with index \ $\beta_0 \in \RR_+$, \ $X_{-1}=0$,
 \ $m_\xi \in \RR_{++}$ \ and \ $m_\eta \in \RR_+$.
\ In case of \ $\beta_0 \in [1, \infty)$, \ assume additionally that there exists
\ $r \in (\beta_0, \infty)$ \ with \ $\EE(\xi^r) < \infty$ \ and \ $\EE(\eta^r) < \infty$.
\ Then for all \ $n \in \NN$,
 \[
   \PP(X_n > x) \sim m_n^{\beta_0} \PP(X_0 > x) \qquad \text{as \ $x \to \infty$,}
 \]
 where \ $m_i$, \ $i \in \ZZ_+$, \ are given in Theorem \ref{Thm3}, and hence, \ $X_n$ \ is also
 regularly varying with index \ $\beta_0$ \ for each \ $n \in \NN$.
\end{Pro}

\noindent{\bf Proof of Proposition \ref{Cor_2GW_X_0_X_-1}.}
Let us fix \ $n \in \NN$.
\ In view of the additive property \eqref{2GW_additive}, it is sufficient to prove
 \[
   \PP\Biggl(\sum_{i=1}^{X_0} \zeta_{i,0}^{(n)} > x\Biggr) \sim m_n^{\beta_0} \PP(X_0 > x) \qquad
   \text{as \ $x \to \infty$.}
 \]
This relation follows from Proposition \ref{FGAMSRS}, since
 \ $\EE(\zeta_{1,0}^{(n)}) =  m_n \in \RR_{++}$, \ $n \in \NN$, \ by \eqref{EXn}.
\proofend

In Appendix \ref{section_2GW}, we present a variant of Proposition \ref{Cor_2GW_X_0_X_-1},
 where the initial values \ $X_{-1}$ \ and \ $X_0$ \ are independent and regularly varying
 together with a second type of proof, see Proposition \ref{2GW_X_0_X_-1}.

\noindent{\bf Proof of Theorem \ref{Thm3}.}
First, note that, by Lemma \ref{help_log_exp}, \ $\EE(\bbone_{\{\vare\ne 0\}} \log(\vare)) < \infty$.
\ We will use the ideas of the proof of Theorem 2.1.1 in Basrak et al.\ \cite{BasKulPal}.
Due to the representation \eqref{2GW_additive}, for each \ $i \in \ZZ_+$, \ we have
 \[
   V_i^{(i)}(\vare_i) \distre \sum_{j=1}^{\vare_i} \zeta_{j,0}^{(i)} ,
 \]
 where \ $\bigl\{\vare_i, \zeta_{j,0}^{(i)} : j \in \NN\bigr\}$ \ are independent random variables
 such that \ $\{\zeta_{j,0}^{(i)} : j \in \NN\}$ \ are independent copies of \ $V_{i,0}$, \ where
 \ $(V_{k,0})_{k\geq-1}$ \ is a second-order Galton--Watson {process (without immigration) with
 initial values \ $V_{0,0} = 1$ \ and \ $V_{-1,0} = 0$, and with the same offspring distributions
 as \ $(X_k)_{k\geq-1}$.
\ For each \ $i \in \ZZ_+$, \ by Proposition \ref{Cor_2GW_X_0_X_-1}, \ we obtain
 \ $\PP(V_i^{(i)}(\vare_i) > x) \sim m_i^\alpha \PP(\vare > x)$ \ as \ $x \to \infty$, \ yielding
 that random variables \ $V_i^{(i)}(\vare_i)$, \ $i \in \ZZ_+$, \ are also regularly varying with
 index \ $\alpha$.
\ Since \ $V_i^{(i)}(\vare_i)$, \ $i \in \ZZ_+$, \ are independent, for each \ $n \in \ZZ_+$, \ by
 Lemma \ref{Lem_conv}, we have
 \begin{align}\label{help2}
  \PP\biggl(\sum_{i=0}^n V_i^{(i)}(\vare_i) > x\biggr)
  \sim \sum_{i=0}^n m_i^\alpha \PP(\vare > x) \qquad \text{as \ $x \to \infty$,}
 \end{align}
 and hence the random variables \ $\sum_{i=0}^n V_i^{(i)}(\vare_i)$, \ $n \in \ZZ_+$, \ are also
 regularly varying with index \ $\alpha$.
\ For each \ $n \in \NN$, \ using that \ $V_i^{(i)}(\vare_i)$, \ $i \in \ZZ_+$, \ are non-negative,
 we have
 \begin{align*}
  \liminf_{x\to\infty} \frac{\pi((x,\infty))}{\PP(\vare>x)}
  &= \liminf_{x\to\infty} \frac{\PP(\sum_{i=0}^\infty V_i^{(i)}(\vare_i)>x)}{\PP(\vare>x)}
   \geq \liminf_{x\to\infty} \frac{\PP(\sum_{i=0}^n V_i^{(i)}(\vare_i)>x)}{\PP(\vare>x)}
   = \sum_{i=0}^n m_i^\alpha ,
 \end{align*}
 hence, letting \ $n \to \infty$, \ we obtain
 \begin{equation}\label{liminf}
  \liminf_{x\to\infty} \frac{\pi((x,\infty))}{\PP(\vare > x)}
  \geq \sum_{i=0}^\infty m_i^\alpha .
 \end{equation}
Moreover, for each \ $n \in \NN$ \ and \ $q \in (0, 1)$, \ we have
 \begin{align*}
  &\limsup_{x\to\infty} \frac{\pi((x,\infty))}{\PP(\vare>x)}
   = \limsup_{x\to\infty}
      \frac{\PP\bigl(\sum_{i=0}^{n-1}V_i^{(i)}(\vare_i)
                     +\sum_{i=n}^\infty V_i^{(i)}(\vare_i)>x\bigr)}
           {\PP(\vare>x)} \\
  &\leq \limsup_{x\to\infty}
         \frac{\PP\bigl(\sum_{i=0}^{n-1}V_i^{(i)}(\vare_i)>(1-q)x\bigr)
   	       +\PP\bigl(\sum_{i=n}^\infty V_i^{(i)}(\vare_i)>qx\bigr)}
              {\PP(\vare>x)}
   \leq L_{1,n}(q) + L_{2,n}(q)
 \end{align*}
 with
 \[
   L_{1,n}(q)
   := \limsup_{x\to\infty}
       \frac{\PP\bigl(\sum_{i=0}^{n-1}V_i^{(i)}(\vare_i)>(1-q)x\bigr)}{\PP(\vare>x)} ,
   \qquad
   L_{2,n}(q)
   := \limsup_{x\to\infty}
       \frac{\PP\bigl(\sum_{i=n}^\infty V_i^{(i)}(\vare_i)>qx\bigr)}{\PP(\vare>x)} .
 \]
Since \ $\vare$ \ is regularly varying with index \ $\alpha$, \ by \eqref{help2}, we obtain
 \begin{gather*}
  L_{1,n}(q)
  = \limsup_{x\to\infty}
     \frac{\PP\bigl(\sum_{i=0}^{n-1}V_i^{(i)}(\vare_i)>(1-q)x\bigr)}{\PP(\vare>(1-q)x)}
     \cdot \frac{\PP(\vare>(1-q)x)}{\PP(\vare>x)}
  = (1 - q)^{-\alpha} \sum_{i=0}^{n-1} m_i^\alpha
 \end{gather*}
 and
 \[
   L_{2,n}(q)
   = \limsup_{x\to\infty}
      \frac{\PP\bigl(\sum_{i=n}^\infty V_i^{(i)}(\vare_i)>qx\bigr)}{\PP(\vare>qx)}
      \cdot \frac{\PP(\vare>qx)}{\PP(\vare>x)}
   = q^{-\alpha}
     \limsup_{x\to\infty}
      \frac{\PP\bigl(\sum_{i=n}^\infty V_i^{(i)}(\vare_i)>qx\bigr)}{\PP(\vare>qx)} ,
 \]
 and hence
 \begin{align*}
  \lim_{n\to\infty} L_{1,n}(q)
  &= (1 - q)^{-\alpha} \sum_{i=0}^\infty m_i^\alpha , \\
  \lim_{n\to\infty} L_{2,n}(q)
  &= q^{-\alpha}
     \lim_{n\to\infty} \limsup_{x\to\infty}
      \frac{\PP\bigl(\sum_{i=n}^\infty V_i^{(i)}(\vare_i)>qx\bigr)}{\PP(\vare>qx)} .
 \end{align*}
The aim of the following discussion is to show
 \begin{equation}\label{limlimsup}
  \lim_{n\to\infty}
   \limsup_{x\to\infty}
    \frac{\PP\bigl(\sum_{i=n}^\infty V_i^{(i)}(\vare_i)>qx\bigr)}{\PP(\vare>qx)}
   = 0, \qquad q \in (0, 1 ).
 \end{equation}

First, we consider the case \ $\alpha \in (0, 1)$.
\ For each \ $x \in \RR_{++}$, \ $n \in \NN$ \ and \ $\delta \in (0, 1)$, \ we have
 \begin{align*}
  &\PP\Biggl(\sum_{i=n}^\infty V_i^{(i)}(\vare_i) > x\Biggr) \\
  &= \PP\Biggl(\sum_{i\geq n} V_i^{(i)}(\vare_i) > x , \;
               \sup_{i\geq n} \varrho^i \vare_i > (1 - \delta) x\Biggr)
     + \PP\Biggl(\sum_{i\geq n} V_i^{(i)}(\vare_i) > x , \;
                 \sup_{i\geq n} \varrho^i \vare_i \leq (1 - \delta) x\Biggr) \\
  &= \PP\Biggl(\sum_{i\geq n} V_i^{(i)}(\vare_i) > x , \;
               \sup_{i\geq n} \varrho^i \vare_i > (1 - \delta) x\Biggr) \\
  &\quad
     + \PP\Biggl(\sum_{i\geq n} V_i^{(i)}(\vare_i)
                 \bbone_{\{\vare_i\leq(1-\delta)\varrho^{-i}x\}} > x , \;
                 \sup_{i\geq n} \varrho^i \vare_i \leq (1 - \delta) x\Biggr) \\
  &\leq \PP\biggl(\sup_{i\geq n} \varrho^i \vare_i > (1 - \delta) x\biggr)
        + \PP\Biggl(\sum_{i\geq n}
                     V_i^{(i)}(\vare_i) \bbone_{\{\vare_i\leq(1-\delta)\varrho^{-i}x\}}
                    > x\Biggr)
   =: P_{1,n}(x, \delta) + P_{2,n}(x, \delta) ,
 \end{align*}
 where \ $\varrho$ \ is given in \eqref{varrho}.
By subadditivity of probability,
 \begin{gather*}
  P_{1,n}(x, \delta)
  \leq \sum_{i\geq n} \PP(\varrho^i \vare_i > (1 - \delta) x)
  = \sum_{i\geq n} \PP(\vare > (1 - \delta) \varrho^{-i} x) .
 \end{gather*}
Using Potter's upper bound (see Lemma \ref{Pb}), for \ $\delta \in(0, \frac{\alpha}{2})$, \ there
 exists \ $x_0 \in \RR_{++}$ \ such that
 \begin{equation}\label{Potter}
  \frac{\PP(\vare>(1-\delta)\varrho^{-i}x)}{\PP(\vare>x)}
  < (1 + \delta) [(1 - \delta) \varrho^{-i}]^{-\alpha + \delta}
   < (1 + \delta) [(1 - \delta) \varrho^{-i}]^{-\frac{\alpha}{2}}
 \end{equation}
 if \ $x \in [x_0, \infty)$ \ and \ $(1 - \delta) \varrho^{-i} \in [1, \infty)$, \ which holds for
 sufficiently large \ $i \in \NN$ \ due to \ $\varrho \in (0, 1)$.
\ Consequently, if \ $\delta \in(0, \frac{\alpha}{2})$, \ then
 \begin{equation*}
  \lim_{n\to\infty} \limsup_{x\to\infty} \frac{P_{1,n}(x,\delta)}{\PP(\vare>x)}
  \leq \lim_{n\to\infty} \sum_{i\geq n} (1 + \delta) [(1 - \delta)\varrho^{-i}]^{-\frac{\alpha}{2}}
  = 0 ,
 \end{equation*}
 since \ $\varrho^{\frac{\alpha}{2}} < 1$ \ (due to \ $\varrho \in (0, 1)$) \ yields
 \ $\sum_{i=0}^\infty (\varrho^{-i})^{-\alpha/2} < \infty$.
\ Now we turn to prove that
 \ $\lim_{n\to\infty} \limsup_{x\to\infty} \frac{P_{2,n}(x,\delta)}{\PP(\vare_1 > x)} = 0$.
\ By Markov's inequality,
 \[
   P_{2,n}(x,\delta)
   \leq \frac{1}{x}
        \sum_{i\geq n}
         \EE\bigl(V_i^{(i)}(\vare_i) \bbone_{\{\vare_i\leq(1-\delta)\varrho^{-i}x\}}\bigr) .
 \]
By the representation \ $V_i^{(i)}(\vare_i) \distre \sum_{j=1}^{\vare_i} \zeta^{(i)}_{j,0}$, \ we
 have
 \begin{align*}
  \EE\bigl(V_i^{(i)}(\vare_i) \bbone_{\{\vare_i\leq(1-\delta)\varrho^{-i}x\}}\bigr)
  &= \EE\Biggl(\sum_{j=1}^{\vare_i}
                \zeta_{j,0}^{(i)} \bbone_{\{\vare_i\leq(1-\delta)\varrho^{-i}x\}}\Biggr)
   = \EE\Biggl[\EE\Biggl(\sum_{j=1}^{\vare_i}
                          \zeta_{j,0}^{(i)}
                          \bbone_{\{\vare_i\leq(1-\delta)\varrho^{-i}x\}}
                         \,\Bigg|\, \vare_i\Biggr)\Biggr] \\
  &= \EE\Biggl(\sum_{j=1}^{\vare_i}
                \EE(\zeta_{1,0}^{(i)}) \bbone_{\{\vare_i\leq(1-\delta)\varrho^{-i}x\}}\Biggr)
   = \EE(\zeta_{1,0}^{(i)})
     \EE\bigl(\vare_i \bbone_{\{\vare_i\leq(1-\delta)\varrho^{-i}x\}}\bigr) ,
 \end{align*}
 since \ $\{\zeta_{j,0}^{(i)} : j \in \NN \}$ \ and \ $\vare_i$ \ are independent.
Moreover,
 \begin{align*}
  \EE\bigl(\vare_i \bbone_{\{\vare_i\leq(1-\delta)\varrho^{-i}x\}}\bigr)
  &= \EE\bigl(\vare \bbone_{\{\vare\leq(1-\delta)\varrho^{-i}x\}}\bigr)
   = \int_0^\infty
      \PP\bigl(\vare \bbone_{\{\vare\leq(1-\delta)\varrho^{-i}x\}} > t\bigr) \, \dd t \\
  &= \int_0^{(1-\delta)\varrho^{-i}x} \PP(t < \vare \leq (1 - \delta) \varrho^{-i} x) \, \dd t
   \leq \int_0^{(1-\delta)\varrho^{-i}x} \PP(\vare > t) \, \dd t .
 \end{align*}
By Karamata's theorem (see, Theorem \ref{Krthm}), we have
 \[
   \lim_{y\to\infty} \frac{\int_0^y\PP(\vare>t)\,\dd t}{y\PP(\vare>y)} = \frac{1}{1-\alpha} ,
 \]
 thus there exists \ $y_0 \in \RR_{++}$ \ such that
 \[
   \int_0^y \PP(\vare > t) \, \dd t \leq \frac{2y\PP(\vare>y)}{1-\alpha} , \qquad
   y \in [y_0, \infty) ,
 \]
 hence
 \[
   \int_0^{(1-\delta)\varrho^{-i}x} \PP(\vare > t) \, \dd t
   \leq \frac{2(1-\delta)\varrho^{-i}x\PP(\vare>(1-\delta)\varrho^{-i}x)}{1-\alpha}
 \]
 whenever \ $(1 - \delta) \varrho^{-i} x \in [y_0, \infty)$, \ which holds for \ $i \geq n$ \ with
 sufficiently large \ $n \in \NN$ \ and \ $x \in [(1-\delta)^{-1} \varrho^n y_0, \infty)$ \ due to
 \ $\varrho \in (0, 1)$.
\ Thus, for sufficiently large \ $n \in \NN$ \ and
 \ $x \in [(1-\delta)^{-1} \varrho^n y_0, \infty)$, \ we obtain
 \begin{align*}
  \frac{P_{2,n}(x,\delta)}{\PP(\vare>x)}
  &\leq \frac{1}{x\PP(\vare>x)}
        \sum_{i\geq n}
        \EE(\zeta_{1,0}^{(i)}) \int_0^{(1-\delta)\varrho^{-i}x} \PP(\vare > t) \, \dd t \\
  &\leq \frac{2(1-\delta)}{1-\alpha}
        \sum_{i\geq n} \frac{\PP(\vare>(1-\delta)\varrho^{-i} x)}{\PP(\vare>x)} ,
 \end{align*}
 since \ $\EE(\zeta_{1,0}^{(i)}) \leq \varrho^i$, \ $i \in \ZZ_+,$ \ by \eqref{1moment_est} and
 \ $\zeta^{(0)}_{1,0} = 1$.
\ Using \eqref{Potter}, we get
 \[
   \frac{P_{2,n}(x,\delta)}{\PP(\vare>x)}
   \leq \frac{2(1-\delta)}{1-\alpha}
        \sum_{i\geq n} (1 + \delta) [(1 - \delta) \varrho^{-i}]^{-\frac{\alpha}{2}}
 \]
 for \ $\delta \in (0, \frac{\alpha}{2})$, \ for sufficiently large \ $n \in \NN$ \ and
 for all \ $x \in [\max(x_0, (1-\delta)^{-1} \varrho^n y_0) , \infty)$.
\ Hence for \ $\delta \in (0, \frac{\alpha}{2})$ \ we have
 \[
   \lim_{n\to\infty} \limsup_{x\to\infty} \frac{P_{2,n}(x,\delta)}{\PP(\vare>x)}
   \leq \lim_{n\to\infty}
         \frac{2(1-\delta^2)}{1-\alpha}
         \sum_{i\geq n} [(1 - \delta) \varrho^{-i}]^{-\frac{\alpha}{2}}
   = 0 ,
 \]
 where the last step follows by the fact that the series
 \ $\sum_{i=0}^\infty (\varrho^i)^{\frac{\alpha}{2}}$ \ is convergent, since
 \ $\varrho \in (0, 1)$.

Consequently, due to the fact that
 \ $\PP(\sum_{i=n}^\infty V_i^{(i)}(\vare_i) > x) \leq P_{1,n}(x,\delta) + P_{2,n}(x,\delta)$,
 \ $x\in\RR_{++}$, \ $n \in \NN$, \ $\delta \in (0, 1)$, \ we obtain \eqref{limlimsup}, and we
 conclude \ $\lim_{n\to\infty} L_{2,n}(q) = 0$ \ for all \ $q \in (0, 1)$.
\ Thus we obtain
 \[
   \limsup_{x\to\infty} \frac{\pi((x,\infty))}{\PP(\vare>x)}
   \leq \lim_{n\to\infty} L_{1,n}(q) + \lim_{n\to\infty} L_{2,n}(q)
   = (1 - q)^{-\alpha} \sum_{i=0}^\infty m_i^\alpha
 \]
 for all \ $q \in (0, 1)$.
\ Letting \ $q \downarrow 0$, \ this yields
 \[
   \limsup_{x\to\infty} \frac{\pi((x,\infty))}{\PP(\vare>x)}
   \leq \sum_{i=0}^\infty m_i^\alpha .
 \]
Taking into account \eqref{liminf}, the proof of \eqref{limlimsup} is complete in case of
 \ $\alpha \in(0, 1)$.

Next, we consider the case \ $\alpha \in [1, 2)$.
\ Note that \eqref{limlimsup} is equivalent to
 \begin{equation*}
  \lim_{n\to\infty}\limsup_{x\to\infty}
   \frac{\PP\bigl(\sum_{i=n}^\infty V_i^{(i)}(\vare_i)>\sqrt{x}\bigr)}{\PP(\vare>\sqrt{x})}
  = \lim_{n\to\infty}\limsup_{x\to\infty}
     \frac{\PP\bigl(\bigl(\sum_{i=n}^\infty V_i^{(i)}(\vare_i)\bigr)^2>x\bigr)}{\PP(\vare^2>x)}
  = 0 .
 \end{equation*}
Repeating a similar argument as for \ $\alpha \in (0, 1)$, \ we obtain
 \begin{align*}
  &\PP\Biggl(\left(\sum_{i=n}^\infty V_i^{(i)}(\vare_i)\right)^2 > x\Biggr) \\
  &= \PP\Biggl(\left(\sum_{i=n}^\infty V_i^{(i)}(\vare_i)\right)^2 > x, \;
                     \sup_{i\geq n} \varrho^{2i} \vare_i^2 > (1 - \delta) x\Biggr)\\
  &\phantom{\quad}
     + \PP\Biggl(\left(\sum_{i=n}^\infty V_i^{(i)}(\vare_i) \right)^2 > x, \;
                       \sup_{i\geq n} \varrho^{2i} \vare_i^2 \leq (1 - \delta) x\Biggr) \\
  &= \PP\Biggl(\left(\sum_{i=n}^\infty V_i^{(i)}(\vare_i)\right)^2 > x, \;
                     \sup_{i\geq n} \varrho^{2i} \vare_i^2 > (1 - \delta) x\Biggr) \\
  &\quad
     + \PP\Biggl(\left(\sum_{i=n}^\infty
                        V_i^{(i)}(\vare_i)
                        \bbone_{\{\vare_i^2\leq(1-\delta)\varrho^{-2i}x\}}\right)^2
                       > x, \;
                       \sup_{i\geq n} \varrho^{2i} \vare_i^2 \leq (1 - \delta) x\Biggr)
 \end{align*}
 \begin{align*}
  &\leq \PP\biggl(\sup_{i\geq n} \varrho^{2i} \vare_i^2 > (1 - \delta) x\biggr)
        + \PP\Biggl(\left(\sum_{i=n}^\infty
                           V_i^{(i)}(\vare_i)
                           \bbone_{\{\vare_i^2\leq(1-\delta)\varrho^{-2i}x\}}\right)^2
                    > x\Biggr) \\
  &=: P_{1,n}(x, \delta) + P_{2,n}(x, \delta)
 \end{align*}
 for each \ $x \in \RR_{++}$, \ $n \in \NN$ \ and \ $\delta \in (0, 1)$.
\ By the subadditivity of probability,
 \[
   P_{1,n}(x, \delta)
   \leq \sum_{i=n}^\infty \PP(\varrho^{2i} \vare_i^2 > (1 - \delta) x)
   = \sum_{i=n}^\infty \PP(\vare^2 > (1 - \delta) \varrho^{-2i} x)
 \]
 for each \ $x \in \RR_{++}$, \ $n \in \NN$ \ and \ $\delta \in (0, 1)$.
\ Since \ $\vare^2$ \ is regularly varying with index \ $\frac{\alpha}{2}$ \ (see Lemma
 \ref{Lem_regl_power}), \ using Potter's upper bound (see Lemma \ref{Pb}) for
 \ $\delta \in \bigl(0, \frac{\alpha}{4}\bigr)$, \ there exists \ $x_0 \in \RR_{++}$ \ such that
 \begin{equation}\label{Potter2}
  \frac{\PP(\vare^2 > (1 - \delta) \varrho^{-2i} x)}{\PP(\vare^2 > x)}
  < (1 + \delta) [(1 - \delta) \varrho^{-2i}]^{-\frac{\alpha}{2} + \delta}
  < (1 + \delta) [(1 - \delta) \varrho^{-2i}]^{-\frac{\alpha}{4}}
 \end{equation}
 if \ $x \in [x_0, \infty)$ \ and \ $(1 - \delta) \varrho^{-2i} \in [1, \infty)$, \ which holds
 for sufficiently large \ $i \in \NN$ \ (due to \ $\varrho \in (0, 1)$).
\ Consequently, if \ $\delta \in (0, \frac{\alpha}{4})$, \ then
 \begin{equation*}
  \lim_{n\to\infty} \limsup_{x\to\infty} \frac{P_{1,n}(x,\delta)}{\PP(\vare^2 >x)}
  \leq \lim_{n\to\infty}
        \sum_{i=n}^\infty (1 + \delta) [(1 - \delta) \varrho^{-2i}]^{-\frac{\alpha}{4}}
  = 0,
 \end{equation*}
 since \ $\varrho^{\frac{\alpha}{2}} < 1$ \ (due to \ $\varrho \in (0, 1)$).
\ By Markov's inequality, for \ $x \in \RR_{++}$, \ $n \in \NN$ \ and \ $\delta \in (0, 1)$, \ we
 have
 \begin{align*}
  \frac{P_{2,n}(x,\delta)}{\PP(\vare^2>x)}
  &\leq \frac{1}{x\PP(\vare^2>x)}
        \EE\Biggl(\left(\sum_{i=n}^\infty V_i^{(i)}(\vare_i)
                         \bbone_{\{\vare_i^2\leq(1-\delta)\varrho^{-2i}x\}}\right)^2\Biggr) \\
  &= \frac{1}{x\PP(\vare^2>x)}
     \EE\Biggl(\sum_{i=n}^\infty V_i^{(i)}(\vare_i)^2
                \bbone_{\{\vare_i^2\leq(1-\delta)\varrho^{-2i}x\}}\Biggr) \\
  &\quad
     + \frac{1}{x\PP(\vare^2>x)}
       \EE\Biggl(\sum_{i,j=n,\; i\ne j}^\infty V_i^{(i)}(\vare_i) V_j^{(j)}(\vare_j)
                  \bbone_{\{\vare_i^2\leq(1-\delta)\varrho^{-2i}x\}}
                  \bbone_{\{\vare_j^2\leq(1-\delta)\varrho^{-2j}x\}}\Biggr) \\
    &=: J_{2,1,n}(x, \delta) + J_{2,2,n}(x, \delta)
 \end{align*}
 for each \ $x \in \RR_{++}$, \ $n \in \NN$ \ and \ $\delta \in (0, 1)$.
By Lemma \ref{Lem_seged_moment}, \eqref{EXn} and \eqref{1moment_est} with \ $X_0 = 1$ \ and
 \ $X_{-1} = 0$, \ we have
 \begin{align*}
  \EE(V_i^{(i)}(n)^2)
  &= \EE\left(\left(\sum_{j=1}^n \zeta_{j,0}^{(i)}\right)^2\right)
   = \sum_{j=1}^n \EE\big((\zeta_{j,0}^{(i)})^2\big)
     + \sum_{j,\ell=1,\,j\ne\ell}^n
        \EE\big(\zeta_{j,0}^{(i)}\big) \EE\big(\zeta_{\ell,0}^{(i)}\big) \\
  &\leq c_\SUB \sum_{j=1}^n \varrho^i + \sum_{j,\ell=1,\,j\ne\ell}^n \varrho^i \varrho^i
   \leq c_\SUB n \varrho^i + (n^2 - n) \varrho^{2i}
   \leq c_\SUB \varrho^i n + \varrho^{2i} n^2
 \end{align*}
 for \ $i, n \in \NN$.
\ Hence, using that
 \ $(\vare_i, V_i^{(i)}(\vare_i))
    \distre \bigl(\vare_i, \sum_{j=1}^{\vare_i} \zeta_{j,0}^{(i)}\bigr)$
 \ and that \ $\vare_i$ \ and \ $\{\zeta^{(i)}_{j,0} : j \in \NN\}$ \ are independent, we have
 \begin{align*}
  J_{2,1,n}(x, \delta)
  &= \sum_{i=n}^\infty
      \frac{\EE\bigl(V_i^{(i)}(\vare_i)^2\bbone_{\{\vare_i^2\leq(1-\delta)\varrho^{-2i}x\}}\bigr)}
           {x\PP(\vare^2>x)} \\
  &= \sum_{i=n}^\infty
      \frac{\EE\bigl(\bigl(\sum_{j=1}^{\vare_i}\zeta_{j,0}^{(i)}\bigr)^2
                     \bbone_{\{\vare_i\leq(1-\delta)^{\frac{1}{2}}\varrho^{-i}x^{\frac{1}{2}}\}}
               \bigr)}
           {x\PP(\vare^2>x)} \\
  &= \sum_{i=n}^\infty
      \frac{\sum_{0\leq\ell\leq(1-\delta)^{\frac{1}{2}}\varrho^{-i}x^{\frac{1}{2}}}
             \EE\Bigl(\bigl(\sum_{j=1}^\ell\zeta_{j,0}^{(i)}\bigr)^2\Bigr)\PP(\vare_i=\ell)}
           {x\PP(\vare^2>x)}\\
  &\leq \sum_{i=n}^\infty
         \frac{\sum_{0\leq\ell\leq(1-\delta)^{\frac{1}{2}}\varrho^{-i}x^{\frac{1}{2}}}
                \left(c_\SUB\varrho^i\ell+\varrho^{2i}\ell^2\right)
                \PP(\vare=\ell)}
              {x\PP(\vare^2>x)}\\
  &= \sum_{i=n}^\infty
      c_\SUB \varrho^i
      \frac{\EE(\vare\bbone_{\{\vare^2\leq(1-\delta)\varrho^{-2i} x\}})}{x\PP(\vare^2>x)}
     + \sum_{i=n}^\infty \varrho^{2i}
        \frac{\EE(\vare^2\bbone_{\{\vare^2\leq(1-\delta)\varrho^{-2i}x\}})}{x\PP(\vare^2>x)} \\
  &=: J_{2,1,1,n}(x, \delta) + J_{2,1,2,n}(x, \delta) .
 \end{align*}
Since \ $\vare^2$ \ is regularly varying with index \ $\frac{\alpha}{2} \in [\frac{1}{2}, 1)$
 \ (see Lemma \ref{Lem_regl_power}), by Karamata's theorem (see, Theorem \ref{Krthm}), we have
 \[
   \lim_{y\to\infty} \frac{\int_0^y \PP(\vare^2>t)\,\dd t}{y\PP(\vare^2>y)}
   = \frac{1}{1-\frac{\alpha}{2}} ,
 \]
 thus there exists \ $y_0 \in \RR_{++}$ \ such that
 \[
   \int_0^y \PP(\vare^2 > t) \, \dd t \leq \frac{2y\PP(\vare^2>y)}{1-\frac{\alpha}{2}} , \qquad
   y \in [y_0, \infty) ,
 \]
 hence
 \begin{align*}
  \EE(\vare^2 \bbone_{\{\vare^2\leq(1-\delta)\varrho^{-2i}x\}})
  &= \int_0^\infty \PP(\vare^2 \bbone_{\{\vare^2\leq(1-\delta)\varrho^{-2i}x\}} > y) \, \dd y\\
  & = \int_0^{(1-\delta)\varrho^{-2i}x}
      \PP(y < \vare^2 \leq (1-\delta) \varrho^{-2i} x) \, \dd y\\
   &\leq \int_0^{(1-\delta)\varrho^{-2i}x} \PP(\vare^2 > t) \, \dd t
     \leq \frac{2(1-\delta)\varrho^{-2i}x\PP(\vare^2>(1-\delta)\varrho^{-2i}x)}
               {1-\frac{\alpha}{2}}
 \end{align*}
 whenever \ $(1 - \delta) \varrho^{-2i} x \in [y_0, \infty)$, \ which holds for \ $i \geq n$
 \ with sufficiently large \ $n \in \NN$, \ and
 \ $x \in [(1 - \delta)^{-1} \varrho^{2n} y_0, \infty)$ \ due to \ $\varrho \in (0, 1)$.
\ Thus for \ $\delta \in (0, \frac{\alpha}{4})$, \ for sufficiently large \ $n \in \NN$
 \ (satisfying \ $(1 - \delta) \varrho^{-2n} \in (1, \infty)$ \ as well) and for all
 \ $x \in [\max(x_0, (1-\delta)^{-1} \varrho^{2n} y_0), \infty)$, \ using \eqref{Potter2}, we
 obtain
 \begin{align*}
  J_{2,1,2,n}(x, \delta)
  &\leq \frac{2(1-\delta)}{1-\frac{\alpha}{2}}
        \sum_{i=n}^\infty \frac{\PP(\vare^2>(1-\delta)\varrho^{-2i}x)}{\PP(\vare^2>x)}
   \leq \frac{2(1-\delta)}{1-\frac{\alpha}{2}}
        \sum_{i=n}^\infty
         (1 + \delta) [(1 - \delta) \varrho^{-2i}]^{-\frac{\alpha}{4}} \\
  &= \frac{2(1-\delta^2)}{1-\frac{\alpha}{2}}
     \sum_{i=n}^\infty [(1 - \delta) \varrho^{-2i}]^{-\frac{\alpha}{4}} .
 \end{align*}
Hence for \ $\delta \in (0, \frac{\alpha}{4})$, \ we have
 \[
   \lim_{n\to\infty} \limsup_{x\to\infty} J_{2,1,2,n}(x, \delta)
   \leq \frac{2(1-\delta^2)}{1-\frac{\alpha}{2}}
        \lim_{n\to\infty} \sum_{i=n}^\infty
         [(1 - \delta) \varrho^{-2i}]^{-\frac{\alpha}{4}}
   = 0 ,
 \]
 yielding \ $\lim_{n\to\infty} \limsup_{x\to\infty} J_{2,1,2,n}(x, \delta) = 0$ \  for
 \ $\delta \in (0, \frac{\alpha}{4})$.
\ Further, if \ $\alpha \in (1, 2)$, \ or \ $\alpha = 1$ \ and \ $m_\vare < \infty$, \ we have
 \[
   J_{2,1,1,n}(x, \delta)
   \leq c_\SUB \sum_{i=n}^\infty \varrho^i \frac{m_\vare}{x\PP(\vare^2>x)} ,
 \]
 and hence, using that \ $\lim_{x\to\infty} x \PP(\vare^2 > x) = \infty$ \ (see Lemma \ref{sv}),
 \begin{align*}
  \lim_{n\to\infty} \limsup_{x\to\infty} J_{2,1,1,n}(x, \delta)
  \leq c_\SUB m_\vare
       \lim_{n\to\infty}
         \Biggl(\sum_{i=n}^\infty \varrho^i \Biggr)
         \limsup_{x\to\infty} \frac{1}{x\PP(\vare^2>x)}
  = 0 ,
 \end{align*}
 yielding \ $\lim_{n\to\infty} \limsup_{x\to\infty} J_{2,1,1,n}(x, \delta) = 0$ \ for
 \ $\delta \in (0, 1)$.

If \ $\alpha = 1$ \ and \ $m_\vare = \infty$, \ then we have
 \[
   J_{2,1,1,n}(x, \delta)
   = \sum_{i=n}^\infty
      c_\SUB \varrho^i
      \frac{\EE\big(\vare
                   \bbone_{\{\vare\leq(1-\delta)^{\frac{1}{2}}\varrho^{-i}x^{\frac{1}{2}}\}} \big)}
           {x\PP(\vare^2>x)}
 \]
 for \ $x \in \RR_{++}$, \ $n \in \NN$ \ and \ $\delta \in (0, 1)$.
\ Note that
 \[
   \EE(\vare \bbone_{\{\vare\leq y\}})
   \leq \int_0^\infty \PP(\vare \bbone_{\{\vare\leq y\}} > t) \, \dd t
   = \int_0^y \PP(t < \vare \leq y) \, \dd t
   \leq \int_0^y \PP(t < \vare) \, \dd t
   =: \tL(y)
 \]
 for \ $y \in \RR_+$.
\ Because of \ $\alpha = 1$, \ Proposition 1.5.9a in Bingham et al.\ \cite{BinGolTeu} yields that
 \ $\tL$ \ is a slowly varying function (at infinity).
By Potter's bounds (see Lemma \ref{Pb}), for every \ $\delta \in \RR_{++}$, \ there exists
 \ $z_0 \in \RR_{++}$ \ such that
 \[
   \frac{\tL(y)}{\tL(z)} < (1 + \delta) \left(\frac{y}{z}\right)^\delta
 \]
 for \ $z \geq z_0$ \ and \ $y \geq z$.
\ Hence, for \ $x \geq z_0^2$, \ we have
 \[
   \EE\bigl(\vare\bbone_{\{\vare\leq(1-\delta)^{\frac{1}{2}}\varrho^{-i}x^{\frac{1}{2}}\}}\bigr)
   \leq \tL\bigl((1 - \delta)^{\frac{1}{2}} \varrho^{-i} x^{\frac{1}{2}}\bigr)
   \leq \tL(\varrho^{-i} x^{\frac{1}{2}})
   \leq (1 + \delta) \varrho^{-i\delta} \widetilde{L}(x^{\frac{1}{2}}) ,
   \qquad i\geq n,
 \]
 where we also used that \ $\tL$ \ is monotone increasing.
Using this, we conclude that for every \ $\delta \in \RR_{++}$, \ there exists \ $z_0 \in \RR_{++}$
 \ such that for \ $x \geq z_0^2$, \ we have
 \[
   J_{2,1,1,n}(x, \delta)
   \leq (1 + \delta) c_\SUB \frac{\tL(x^{\frac{1}{2}})}{x\PP(\vare^2>x)}
        \sum_{i=n}^\infty \varrho^{-i\delta} .
 \]
Here, since \ $\varrho \in (0, 1)$ \ and \ $\delta \in \RR_{++}$, \ we have
 \ $\lim_{n\to\infty} \sum_{i=n}^\infty \varrho^{-i\delta} = 0$, \ and
 \[
  \frac{\tL(\sqrt{x})}{x\PP(\vare^2>x)}
  = \frac{\tL(\sqrt{x})}{x^{1/4}} \cdot \frac{1}{x^{3/4} \PP(\vare > \sqrt{x})}
  \to 0 \qquad \text{as \ $x \to \infty$,}
 \]
 by Lemma \ref{sv}, due to the fact that \ $\tL$ \ is slowly varying and the function
 \ $\RR_{++} \ni x \mapsto \PP(\vare > \sqrt{x})$ \ is regularly varying with index \ $-1/2$.
\ Hence \ $\lim_{n\to\infty} \limsup_{x\to\infty} J_{2,1,1,n}(x, \delta) = 0$ \ for
 \ $\delta \in (0, 1)$ \ in case of \ $\alpha = 1$ \ and \ $m_\vare = \infty$.

Consequently, \ we have \ $\lim_{n\to\infty} \limsup_{x\to\infty} J_{2,1,n}(x, \delta) = 0$ \ for
 \ $\delta \in (0, \frac{\alpha}{4})$.

Now we turn to prove \ $\lim_{n\to\infty} \limsup_{x\to\infty} J_{2,2,n}(x, \delta) = 0$ \ for
 \ $\delta \in (0, 1)$.
\ Using that \ $\{(\vare_i, V_i^{(i)}(\vare_i) : i \in \NN\}$ \ are independent, we have
 \begin{align*}
  J_{2,2,n}(x, \delta)
  \leq \frac{1}{x\PP(\vare^2>x)}
       \sum_{i,j=n,\,i\ne j}^\infty
        \EE(V_i^{(i)}(\vare_i)
        \bbone_{\{\vare_i^2\leq(1-\delta)\varrho^{-2i}x\}})
        \EE(V_j^{(j)}(\vare_j) \bbone_{\{\vare_j^2\leq(1-\delta)\varrho^{-2j}x\}}) .
 \end{align*}
Here, using that
 \ $\left(\vare_i, V_i^{(i)}(\vare_i)\right)
    \distre \bigl(\vare_i, \sum_{j=1}^{\vare_i} \zeta_{j,0}^{(i)}\bigr)$,
 \ where \ $\vare_i$ \ and \ $\{\zeta^{(i)}_{j,0} : j \in \NN\}$ \ are independent, and
 \eqref{1moment_est} with \ $X_0 = 1$ \ and \ $X_{-1} = 0$, \ we have
 \begin{align*}
  \EE(V_i^{(i)}(\vare_i) \bbone_{\{\vare_i^2\leq(1-\delta)\varrho^{-2i}x\}})
  &= \EE\left(\left(\sum_{j=1}^{\vare_i} \zeta_{j,0}^{(i)}\right)
              \bbone_{\{\vare_i^2\leq(1-\delta)\varrho^{-2i}x\}}\right) \\
  &= \sum_{\ell=0}^{\lfloor(1-\delta)^{\frac{1}{2}}\varrho^{-i}x^{\frac{1}{2}}\rfloor}
      \EE\left(\sum_{j=1}^{\ell} \zeta_{j,0}^{(i)}\right) \PP(\vare_i = \ell) \\
  &\le \sum_{\ell=0}^{\lfloor(1-\delta)^{\frac{1}{2}}\varrho^{-i}x^{\frac{1}{2}}\rfloor}
        \ell \varrho^i \PP(\vare_i = \ell)
   = \varrho^i \EE(\vare_i \bbone_{\{\vare_i^2\leq(1-\delta)\varrho^{-2i}x\}})
 \end{align*}
 for \ $x \in \RR_{++}$ \ and \ $\delta \in (0, 1)$.
\ If \ $\alpha \in (1, 2)$, \ or \ $\alpha = 1$ \ and \ $m_\vare < \infty$, \ then
 \begin{align*}
  J_{2,2,n}(x, \delta)
  &\leq \frac{1}{x\PP(\vare^2>x)}
        \sum_{i,j=n,\,i\ne j}^\infty
         \varrho^{i+j}
         \EE(\vare_i \bbone_{\{\vare_i^2\leq(1-\delta)\varrho^{-2i}x\}})
         \EE(\vare_j \bbone_{\{\vare_j^2\leq(1-\delta)\varrho^{-2j}x\}}) \\
  &\leq \frac{m_\vare^2}{x\PP(\vare^2>x)}
        \sum_{i,j=n,\,i\ne j}^\infty \varrho^{i+j}
   \leq \frac{m_\vare^2}{x\PP(\vare^2>x)}
        \left(\sum_{i=n}^\infty \varrho^i\right)^2
 \end{align*}
 for \ $x \in \RR_{++}$ \ and \ $\delta \in (0, 1)$, \ and then, by  Lemma \ref{sv},
 \begin{align*}
  \lim_{n\to\infty} \limsup_{x\to\infty} J_{2,2,n}(x, \delta)
  &\leq m_\vare^2
        \lim_{n\to\infty}
         \left(\sum_{i=n}^\infty \varrho^i\right)^2
         \limsup_{x\to\infty} \frac{1}{x\PP(\vare^2>x)} \\
  &= m_\vare^2 \left(\lim_{n\to\infty} \frac{\varrho^{2n}}{(1-\varrho)^2}\right) \cdot 0
   = 0 ,
 \end{align*}
 yielding that \ $\lim_{n\to\infty} \limsup_{x\to\infty} J_{2,2,n}(x, \delta) = 0$.

If \ $\alpha = 1$ \ and \ $m_\vare = \infty$, \ then we can apply the same argument as for
 \ $J_{2,1,1,n}(x, \delta)$.
\ Namely,
 \begin{align*}
  J_{2,2,n}(x, \delta)
  &\leq \frac{(1+\delta)^2}{x\PP(\vare^2>x)}
        \sum_{i,j= n,\,i\ne j}^\infty
         \varrho^{(1-\delta)(i+j)} (\widetilde L(x^{\frac{1}{2}}))^2 \\
  &\leq (1+\delta)^2 \frac{(\tL(x^{\frac{1}{2}}))^2}{x\PP(\vare^2>x)}
        \sum_{i,j=n,\,i\ne j}^\infty \varrho^{(1-\delta)(i+j)}
   = (1+\delta)^2 \frac{(\tL(x^{\frac{1}{2}}))^2}{x\PP(\vare^2>x)}
     \left(\sum_{i=n}^\infty \varrho^{(1-\delta)i}\right)^2
 \end{align*}
 for \ $x \in \RR_{++}$ \ and \ $\delta \in (0, 1)$, \ where
 \[
   \frac{(\tL(x^{\frac{1}{2}}))^2}{x\PP(\vare^2 >x)}
   = \left(\frac{\tL(x^{\frac{1}{2}})}{x^{\frac{1}{2}}}\right)^2
     \frac{1}{x^{\frac{3}{4}} \PP(\vare > \sqrt{x})}
   \to 0 \qquad \text{as \ $x \to \infty$,}
 \]
 yielding that \ $\lim_{n\to\infty} \limsup_{x\to\infty} J_{2,2,n}(x, \delta) = 0$ \ for
 \ $\delta \in (0, 1)$ \ in case of \ $\alpha = 1$ \ and \ $m_\vare = \infty$ \ as well.

Consequently,
 \ $\lim_{n\to\infty} \limsup_{x\to\infty} \frac{P_{2,n}(x,\delta)}{\PP(\vare^2>x)} = 0$ \ for
 \ $\delta \in (0, \frac{\alpha}{4})$ \ yielding \eqref{limlimsup} in case of \ $\alpha \in [1, 2)$
 \ as well, and we conclude \ $\lim_{n\to\infty} L_{2,n}(q) = 0$ \ for all \ $q \in (0, 1)$.
\ The proof can be finished as in case of \ $\alpha \in (0, 1)$.
\proofend

\begin{Rem}
The statement of Theorem \ref{Thm3} remains true in the case when \ $m_\xi \in (0, 1)$ \ and
 \ $m_\eta = 0$.
\ In this case we get the statement for classical Galton--Watson processes, see Theorem 2.1.1 in
 Basrak et al.\ \cite{BasKulPal} or Theorem \ref{GWI_stat}.
However, note that this is not a special case of Theorem \ref{Thm3}, since in this case the mean
 matrix \ $\bM_{\xi,\eta}$ \ is not primitive.
\proofend
\end{Rem}

\appendix

\vspace*{5mm}

\noindent{\bf\Large Appendices}

\section{Representations of second-order Galton--Watson processes without or with immigration}
\label{section_prel}

First, we recall a representation of a second-order Galton--Watson process without or with
 immigration as a (special) 2-type Galton--Watson process without or with immigration,
 respectively.
Let \ $(X_n)_{n\geq -1}$ \ be a second-order Galton--Watson process with immigration given in
 \eqref{2oGWI}, and let us introduce the random vectors
 \begin{align}\label{help15}
  \bY_n := \begin{bmatrix}
            Y_{n,1} \\
            Y_{n,2} \\
           \end{bmatrix}
        := \begin{bmatrix}
            X_n \\
            X_{n-1} \\
           \end{bmatrix} ,
  \qquad n \in \ZZ_+ .
 \end{align}
Then we have
 \begin{align}\label{help16}
  \bY_n = \sum_{i=1}^{Y_{n-1,1}}
           \begin{bmatrix} \xi_{n,i} \\ 1 \end{bmatrix}
          + \sum_{j=1}^{Y_{n-1,2}}
             \begin{bmatrix} \eta_{n,j} \\ 0 \end{bmatrix}
          + \begin{bmatrix} \vare_n \\ 0 \end{bmatrix} ,
  \qquad n \in \NN ,
 \end{align}
 hence \ $(\bY_n)_{n\in\ZZ_+}$ \ is a (special) 2-type Galton--Watson process with immigration and
 with initial vector
 \[
   \bY_0 = \begin{bmatrix} X_0 \\ X_{-1} \end{bmatrix} .
 \]
In fact, the type  1 and  2 individuals are identified with individuals of age \ $0$ \ and \ $1$,
 \ respectively, and for each \ $n, i, j \in \NN$, \ at time \ $n$, \ the \ $i^\mathrm{th}$
 \ individual of type  1 of the \ $(n-1)^\mathrm{th}$ \ generation produces \ $\xi_{n,i}$
 \ individuals of type  1 and exactly one individual of type 2, and the \ $j^\mathrm{th}$
 \ individual of type 2 of the \ $(n-1)^\mathrm{th}$ \ generation produces \ $\eta_{n,j}$
 \ individuals of type 1 and no individual of type 2.

The representation \eqref{help16} works backwards as well, namely, let \ $(\bY_k)_{k \in \ZZ_+}$
 \ be a special 2-type Galton--Watson process with immigration  given by
 \begin{align}\label{help17}
  \bY_k = \sum_{j=1}^{Y_{k-1,1}}
           \begin{bmatrix} \xi_{k,j,1,1} \\ 1 \end{bmatrix}
          + \sum_{j=1}^{Y_{k-1,2}}
             \begin{bmatrix} \xi_{k,j,2,1} \\ 0 \end{bmatrix}
          + \begin{bmatrix} \vare_{k,1} \\ 0 \end{bmatrix} ,
  \qquad k \in \NN ,
 \end{align}
 where \ $\bY_0$ \ is a \ 2-dimensional integer-valued random vector.
Here, for each \ $k, j \in \NN$ \ and \ $i \in \{1,2\}$, \ $\xi_{k,j,i,1}$ \ denotes the number of
 type 1 offsprings in the \ $k^\mathrm{th}$ generation produced by the \ $j^\mathrm{th}$ offspring
 of the \ $(k-1)^\mathrm{th}$ generation of type \ $i$, \ and \ $\vare_k$ \ denotes the number of
 type 1 immigrants in the \ $k^\mathrm{th}$ generation.
For the second coordinate process of \ $(\bY_k)_{k\in\ZZ_+}$, \ we get \ $Y_{k,2} = Y_{k-1,1}$,
 \ $k \in \NN$, \ and substituting this into \eqref{help17}, the first coordinate process of
 \ $(\bY_k)_{ k\in\ZZ_+}$ \ satisfies
 \[
   Y_{k,1} = \sum_{j=1}^{Y_{k-1,1}} \xi_{k,j,1,1} + \sum_{j=1}^{Y_{k-2,1}} \xi_{k,j,2,1}
             + \vare_{k,1}, \qquad  k\geq 2.
 \]
Thus, the first coordinate process of \ $(\bY_k)_{k\in\ZZ_+}$ \ given by \eqref{help17} satisfies
 equation \eqref{2oGWI} with \ $X_n := Y_{n,1}$, \ $\xi_{n,i}:=\xi_{n,i,1,1}$, \ $\eta_{n,j}:= \xi_{n,j,2,1}$, \ $\vare_n := \vare_{n,1}$,
 \ $n, i, j \in \NN$, \ and with initial values \ $X_0:= Y_{0,1}$ \ and \ $X_{-1} := Y_{0,2}$, \ i.e.,
 it is a second-order Galton--Watson process with immigration.
Moreover, the second coordinate process of \ $(\bY_k)_{k\in\ZZ_+}$ \ also satisfies equation
 \eqref{2oGWI} with \ $X_n := Y_{n+1,2}$, \ $\xi_{n,i}:=\xi_{n,i,1,1}$, \ $\eta_{n,j}:= \xi_{n,j,2,1}$, \ $\vare_n := \vare_{n,1}$,
 \ $n, i, j \in \NN$, \ and with initial values \ $X_0:= Y_{0,1}$
 \ and \ $X_{-1} := Y_{0,2}$, \ i.e., it is also a second-order Galton--Watson process with
 immigration.

Note that, for a second-order Galton--Watson process \ $(X_n)_{n\geq -1}$ \ (without immigration),
 the additive (or branching) property of a 2-type Galton--Watson process (without immigration)
 (see, e.g. in Athreya and Ney \cite[Chapter V, Section 1]{AthNey}), together with the law of total
 probability, for each \ $n \in \NN$, \ imply
 \begin{equation}\label{2GW_additive}
  X_n \distre \sum_{i=1}^{X_0} \zeta_{i,0}^{(n)} + \sum_{j=1}^{X_{-1}} \zeta_{j,-1}^{(n)} ,
 \end{equation}
 where \ $\bigl\{ (X_0, X_{-1}), \zeta_{i,0}^{(n)}, \zeta_{j,-1}^{(n)} : i, j \in \NN\bigr\}$ \ are
 independent random variables such that \ $\{\zeta_{i,0}^{(n)} : i \in \NN\}$ \ are independent
 copies of \ $V_{n,0}$ \ and \ $\{\zeta_{j,-1}^{(n)} : j \in \NN\}$ \ are independent copies of
 \ $V_{n,-1}$, \ where \ $(V_{k,0})_{k\geq-1}$ \ and \ $(V_{k,-1})_{k\geq-1}$ \ are second-order
 Galton--Watson processes (without immigration) with initial values \ $V_{0,0} = 1$,
 \ $V_{-1,0} = 0$, \ $V_{0,-1} = 0$ \ and \ $V_{-1,-1} = 1$, \ and with the same offspring
 distributions as \ $(X_k)_{k\geq-1}$.

Moreover, if \ $(X_n)_{n\geq -1}$ \ is a second-order Galton--Watson process with immigration, then
 for each \ $n \in \NN$, \ we have
 \begin{equation}\label{2GWI_additive}
  X_n = V_0^{(n)}(X_0, X_{-1}) + \sum_{i=1}^n V_i^{(n-i)}(\vare_i, 0) ,
 \end{equation}
 where \ $\bigl\{V_0^{(n)}(X_0, X_{-1}), V_i^{(n-i)}(\vare_i, 0) : i \in \{1, \ldots, n\}\bigr\}$
 \ are independent random variables such that \ $V_0^{(n)}(X_0, X_{-1})$ \ represents the number of
 newborns at time \ $n$, \ resulting from the initial individuals \ $X_0$ \ at time \ $0$ \ and
 \ $X_{-1}$ \ at time \ $-1$, \ and for each \ $i \in \{1, \ldots, n\}$,
 \ $V_i^{(n-i)}(\vare_i, 0)$ \ represents the number of newborns at time \ $n$, \ resulting from
 the immigration \ $\vare_i$ \ at time \ $i$.
\ Indeed, considering the (special) 2-type Galton--Watson process \ $(\bY_k)_{k\in\ZZ_+}$ \ with
 immigration given in \eqref{help15} and applying formula (1.1) in Kaplan \cite{Kap2}, we obtain
 \begin{equation}\label{2typeGWI_additive}
  \bY_n = \bV_0^{(n)}(\bY_0) + \sum_{i=1}^n \bV_i^{(n-i)}(\bvare_i) \qquad \text{with} \qquad
  \bvare_i := \begin{bmatrix} \vare_i \\ 0 \end{bmatrix} , \quad i \in \NN ,
 \end{equation}
 for each \ $n \in \NN$, \ where
 \ $\bigl\{\bV_0^{(n)}(\bY_0), \bV_i^{(n-i)}(\bvare_i) : i \in \{1, \ldots, n\}\bigr\}$ \ are
 independent random vectors such that \ $\bV_0^{(n)}(\bY_0)$ \ represents the number of individuals
 alive at time \ $n$, \ resulting from the initial individuals \ $\bY_0$ \ at time \ $0$, \ and for
 each \ $i \in \{1, \ldots, n\}$, \ $\bV_i^{(n-i)}(\bvare_i)$ \ represents the number of
 individuals alive at time \ $n$, \ resulting from the immigration \ $\bvare_i$ \ at time \ $i$.
\ Clearly, \ $(\bV_0^{(k)}(\bY_0))_{k\in\ZZ_+}$ \ and \ $(\bV_i^{(k)}(\bvare_i))_{k\in\ZZ_+}$,
 \ $i \in \{1, \ldots, n\}$, \ are (special) 2-type Galton--Watson processes (without immigration)
 of the form \eqref{help17} with initial vectors \ $\bV_0^{(0)}(\bY_0) = \bY_0$ \ and
 \ $\bV_i^{(0)}(\bvare_i) = \bvare_i$, \ $i \in \{1, \ldots, n\}$, \ respectively, and with the
 same offspring distributions as \ $(\bY_k)_{k\in\ZZ_+}$.
\ Using the considerations for the backward representation presented before, the first coordinates in
 \eqref{2typeGWI_additive} gives \eqref{2GWI_additive}, where
 \ $(V_0^{(k)}(X_0, X_{-1}))_{k\geq-1}$ \ and \ $(V_i^{(k)}(\vare_i, 0))_{k\geq-1}$,
 \ $i \in \{1, \ldots, n\}$, \ are second-order Galton--Watson processes (without immigration) with
 initial values \ $V_0^{(0)}(X_0, X_{-1}) = X_0$, \ $V_0^{(-1)}(X_0, X_{-1}) = X_{-1}$,
 \ $V_i^{(0)}(\vare_i, 0) = \vare_i$ \ and \ $V_i^{(-1)}(\vare_i, 0) = 0$,
 \ $i \in \{1, \ldots, n\}$, \ and with the same offspring distributions as \ $(X_k)_{k\geq-1}$.

\section{On the expectation of second-order Galton--Watson processes with immigration}
\label{App0}

Our aim is to derive an explicit formula for the expectation of a second-order Galton--Watson
 process with immigration at time \ $n$ \ and to describe its asymptotic behavior as
 \ $n \to \infty$.

Recall that \ $\xi$, \ $\eta$ \ and \ $\vare$ \ are random variables such that
 \ $\xi \distre \xi_{1,1}$, \ $\eta \distre \eta_{1,1}$ \ and \ $\vare \distre \vare_1$, \ and we
 put \ $m_\xi = \EE(\xi) \in [0, \infty]$, \ $m_\eta = \EE(\eta) \in [0, \infty]$ \ and
 \ $m_\vare = \EE(\vare) \in [0, \infty]$.
\ If \ $m_\xi \in \RR_+$, \ $m_\eta \in \RR_+$, \ $m_\vare \in \RR_+$, \ $\EE(X_0) \in \RR_+$ \ and
 \ $\EE(X_{-1}) \in \RR_+$, \ then \eqref{2oGWI} implies
 \begin{equation*}
  \EE(X_n \mid \cF_{n-1}^X) = X_{n-1} m_\xi + X_{n-2} m_\eta + m_\vare , \qquad  n \in \NN ,
 \end{equation*}
 where \ $\cF_n^X := \sigma(X_{-1}, X_0, \dots, X_n)$, \ $n \in \ZZ_+$.
\ Consequently,
 \[
   \EE(X_n) = m_\xi \EE(X_{n-1}) + m_\eta \EE(X_{n-2}) + m_\vare , \qquad n \in \NN ,
 \]
 which can be written in the matrix form
 \begin{equation}\label{recEX_n}
  \begin{bmatrix}
   \EE(X_n) \\
   \EE(X_{n-1})
  \end{bmatrix}
  = \bM_{\xi,\eta}
    \begin{bmatrix}
     \EE(X_{n-1}) \\
     \EE(X_{n-2})
    \end{bmatrix}
    + \begin{bmatrix}
       m_\vare \\
       0
      \end{bmatrix} , \qquad n \in \NN ,
 \end{equation}
 with
 \begin{equation}\label{bA}
  \bM_{\xi,\eta} := \begin{bmatrix}
          m_\xi & m_\eta \\
          1 & 0
         \end{bmatrix} .
 \end{equation}
Note that \ $\bM_{\xi,\eta}$ \ is the mean matrix of the 2-type Galton--Watson process
 \ $(\bY_n)_{n\in\ZZ_+}$ \ given in \eqref{help15}.
Thus, we conclude
 \begin{equation}\label{EX}
  \begin{bmatrix}
   \EE(X_n) \\
   \EE(X_{n-1})
  \end{bmatrix}
  = \bM_{\xi,\eta}^n
    \begin{bmatrix}
     \EE(X_0) \\
      \EE(X_{-1})
    \end{bmatrix}
    + \sum_{k=1}^n
       \bM_{\xi,\eta}^{n-k}
       \begin{bmatrix}
        m_\vare \\
        0
       \end{bmatrix} ,
  \qquad n \in \NN .
 \end{equation}
Hence, the asymptotic behavior of the sequence \ $(\EE(X_n))_{n\in\NN}$ \ depends on the asymptotic
 behavior of the powers \ $(\bM_{\xi,\eta}^n)_{n\in\NN}$, \ which is related to the spectral
 radius \ $\varrho$ \ of \ $\bM_{\xi,\eta}$, \ see Lemma \ref{asymptotics_EX_n} and
 \eqref{varrho}.
If \ $(X_n)_{n\geq-1}$ \ is a second-order Galton--Watson process with immigration such that
 \ $m_\xi \in \RR_+$ \ and \ $m_\eta \in \RR_+$, \ then \ $(X_n)_{n\geq-1}$ \ is called
 subcritical, critical or supercritical if \ $\varrho < 1$, \ $\varrho = 1$ \ or \ $\varrho > 1$,
 \ respectively.
It is easy to check that a second-order Galton--Watson process with immigration is subcritical,
 critical or supercritical if and only if \ $m_\xi + m_\eta < 1$, \ $m_\xi + m_\eta = 1$ \ or
 \ $m_\xi + m_\eta > 1$, \ respectively.
We call the attention that for the classification of second-order Galton--Watson process with
 immigration we do not suppose the finiteness of the expectation of \ $X_0$, \ $X_{-1}$ \ or
 \ $\vare$.

\begin{Lem}\label{asymptotics_EX_n}
Let \ $(X_n)_{n\geq-1}$ \ be a second-order Galton--Watson process with immigration such that
 \ $m_\xi \in \RR_+$, \ $m_\eta \in \RR_+$, \ $m_\vare \in \RR_+$, \ $\EE(X_0) \in \RR_+$
 \ and \ $\EE(X_{-1}) \in \RR_+$.

If \ $m_\xi = 0$ \ and \ $m_\eta = 0$, \ then, for all \ $n \in \NN$, \ we have
 \ $\EE(X_n) = m_\vare$.

If \ $m_\xi + m_\eta > 0$, \ then, for all \ $n \in \NN$, \ we have
 \begin{equation}\label{EXn}
  \EE(X_n) = \frac{\lambda_+^{n+1}-\lambda_-^{n+1}}{\lambda_+-\lambda_-} \EE(X_0)
             + \frac{\lambda_+^n-\lambda_-^n}{\lambda_+-\lambda_-} m_\eta \EE(X_{-1})
             + \frac{C_n(\lambda_+,\lambda_-)}{\lambda_+-\lambda_-} m_\vare
 \end{equation}
 with
 \[
   C_n(\lambda_+,\lambda_-)
   := \begin{cases}
       \lambda_+\frac{1-\lambda_+^n}{1-\lambda_+} - \lambda_- \frac{1-\lambda_-^n}{1-\lambda_-}
        & \text{if \ $\lambda_+ \ne 1$,} \\[1mm]
       n - \lambda_- \frac{1-\lambda_-^n}{1-\lambda_-} & \text{if \ $\lambda_+ = 1$,}
      \end{cases}
 \]
 where \ $\lambda_+$ \ and \ $\lambda_-$ \ are given in \eqref{m_n}, and hence
 \begin{align*}
  \EE(X_n)
  = \begin{cases}
     \frac{m_\vare}{(1-\lambda_+)(1-\lambda_-)} + \OO(\lambda_+^n)
      & \text{if \ $\lambda_+ \in (0, 1)$,} \\[1mm]
     \frac{m_\vare}{1-\lambda_-} n + \OO(1) & \text{if \ $\lambda_+ = 1$,} \\[1mm]
     \frac{1}{\lambda_+-\lambda_-}
     \bigl(\lambda_+ \EE(X_0) + m_\eta \EE(X_{-1}) + \frac{\lambda_+}{\lambda_+-1}m_\vare \bigr)
     \lambda_+^n
     + \OO(1+|\lambda_-|^n)
      & \text{if \ $\lambda_+ \in (1, \infty)$}
    \end{cases}
 \end{align*}
 as \ $n \to \infty$.
\ Moreover, \ $\lambda_+$ \ is the spectral radius \ $\varrho$ \ of \ $\bM_{\xi,\eta}$.

Further, in case of \ $m_\vare = 0$, \ we have the following more precise statements:

If \ $m_\xi = 0$, \ $m_\eta > 0$ \ and \ $m_\vare = 0$, \ then, for all \ $k \in \NN$, \ we have
 \ $\EE(X_{2k-1}) = \EE(X_{-1}) \lambda_+^{2k}$ \ and \ $\EE(X_{2k}) = \EE(X_0) \lambda_+^{2k}$.

If \ $m_\xi > 0$, \ $m_\eta = 0$ \ and \ $m_\vare = 0$, \ then, for all \ $n \in \NN$, \ we have
 \ $\EE(X_n) = \EE(X_0) \lambda_+^n$.

If \ $m_\xi > 0$, \ $m_\eta > 0$ \ and \ $m_\vare = 0$, \ then
 \[
   \EE(X_n) = \frac{\lambda_+\EE(X_0)+m_\eta\EE(X_{-1})}{\lambda_+-\lambda_-} \lambda_+^n
              + \OO(|\lambda_-|^n)
   \qquad \text{as \ $n \to \infty$.}
 \]

If \ $m_\vare = 0$, \ i.e., there is no immigration, then
 \begin{equation}\label{1moment_est}
  \EE(X_n) \leq \varrho^n \EE(X_0)  + \varrho^{n-1} m_\eta \EE(X_{-1}), \qquad n \in \NN .
 \end{equation}
\end{Lem}

\noindent{\bf Proof.}
We are going to use \eqref{EX}.
The matrix \ $\bM_{\xi,\eta}$ \ has eigenvalues
 \[
   \lambda_+ = \frac{m_\xi + \sqrt{m_\xi^2 + 4 m_\eta}}{2} , \qquad
   \lambda_- = \frac{m_\xi - \sqrt{m_\xi^2 + 4 m_\eta}}{2} ,
 \]
 satisfying \ $\lambda_+ \in \RR_+$ \ and \ $\lambda_- \in [-\lambda_+, 0]$, \ hence the spectral
 radius of \ $\bM_{\xi,\eta}$ \ is
 \begin{equation}\label{varrho}
  \varrho = \lambda_+ = \frac{m_\xi + \sqrt{m_\xi^2 + 4m_\eta}}{2}.
 \end{equation}
In what follows, we suppose that \ $m_\xi + m_\eta > 0$, \ which yields that
 \ $\lambda_+ \in \RR_{++}$ \ and \ $\lambda_- \in (-\lambda_+, 0]$.
\ One can easily check that the powers of \ $\bM_{\xi,\eta}$ \ can be written in the form
 \begin{equation}\label{powerM}
  \bM_{\xi,\eta}^n
  = \frac{\lambda_+^n}{\lambda_+-\lambda_-}
    \begin{bmatrix}
     \lambda_+ & m_\eta \\
     1 & -\lambda_-
    \end{bmatrix}
    + \frac{\lambda_-^n}{\lambda_+-\lambda_-}
      \begin{bmatrix}
       -\lambda_- & -m_\eta \\
       -1 & \lambda_+
      \end{bmatrix} ,
  \qquad n \in  \ZZ_+.
 \end{equation}
Consequently,
 \[
   \varrho^{-n} \bM_{\xi,\eta}^n
   \to  \frac{1}{\lambda_+-\lambda_-}\begin{bmatrix}
        \lambda_+ & m_\eta \\
        1 & -\lambda_-
       \end{bmatrix}
   \qquad \text{as \ $n \to \infty$.}
 \]
Moreover, \eqref{EX} and \eqref{powerM} yield
 \begin{align*}
  \begin{bmatrix}
   \EE(X_n) \\
   \EE(X_{n-1})
  \end{bmatrix}
  &= \frac{\EE(X_0)}{\lambda_+-\lambda_-}
     \begin{bmatrix}
      \lambda_+^{n+1} - \lambda_-^{n+1} \\
      \lambda_+^n - \lambda_-^n
     \end{bmatrix}
    + \frac{\EE(X_{-1})}{\lambda_+ - \lambda_{-}}
      \begin{bmatrix}
       m_\eta(\lambda_+^n - \lambda_-^n) \\
       -\lambda_-\lambda_+^n + \lambda_+\lambda_-^n
      \end{bmatrix} \\
  &\quad
    + \frac{m_\vare}{\lambda_+-\lambda_-}
      \sum_{k=1}^n
       \begin{bmatrix}
        \lambda_+^{n-k+1} - \lambda_-^{n-k+1} \\
        \lambda_+^{n-k} - \lambda_-^{n-k}
       \end{bmatrix} , \qquad n \in \ZZ_+ ,
 \end{align*}
 and hence, we obtain \eqref{EXn} and \eqref{1moment_est}.
Indeed, by \eqref{powerM} and by \ $\lambda_+  \in\RR_{++}$ \ and
 \ $-\lambda_+ < \lambda_- \leq 0$, \ for each \ $k \in \ZZ_+$, \ we have
 \[
   \frac{\lambda_+^{2k+1}-\lambda_-^{2k+1}}{\lambda_+-\lambda_-}
   = \sum_{i=0}^{2k} \lambda_-^i \lambda_+^{2k-i}
   = \lambda_+^{2k}
     + \sum_{j=1}^{k}
        (\lambda_-^{2j-1} \lambda_+^{2k-2j+1} + \lambda_-^{2j} \lambda_+^{2k-2j})
   \leq \lambda_+^{2k} ,
 \]
 since
 \ $\lambda_-^{2j-1} \lambda_+^{2k-2j+1} + \lambda_-^{2j} \lambda_+^{2k-2j}
    = \lambda_-^{2j-1} \lambda_+^{2k-2j} (\lambda_+ + \lambda_-) \leq 0$,
 \ and, in a similar way,
 \[
   \frac{\lambda_+^{2k+2}-\lambda_-^{2k+2}}{\lambda_+-\lambda_-}
   = \sum_{i=0}^{2k+1} \lambda_-^i \lambda_+^{2k+1-i}
   = \lambda_+^{2k+1}
     + \sum_{j=1}^{k}
        (\lambda_-^{2j-1} \lambda_+^{2k-2j+1} + \lambda_-^{2j} \lambda_+^{2k-2j})
     + \lambda_-^{2k+1}
   \leq \lambda_+^{2k+1} .
 \]
Further, if \ $\lambda_+\in(0,1)$, \ then
 \begin{align*}
   \frac{m_\vare}{\lambda_+ - \lambda_-} C_n(\lambda_+, \lambda_-)
   &= \frac{m_\vare}{\lambda_+ - \lambda_-} \frac{\lambda_+ - \lambda_- + \lambda_+^{n+1}(\lambda_- - 1) + \lambda_-^{n+1}(1-\lambda_+)}
                                              {(1-\lambda_+)(1-\lambda_-)} \\
   &= \frac{m_\vare}{(1-\lambda_+)(1-\lambda_-)} + \OO(\lambda_+^n)
 \end{align*}
 as \ $n \to \infty$.
\ The other statements easily follow from \eqref{EXn}.
\proofend

\section{Moment estimations}
\label{App_moments}

The first moment of a second-order Galton--Watson process \ $(X_n)_{ n\geq -1}$ \ (without
 immigration) can be estimated by \eqref{1moment_est}.
Next, we present an auxiliary lemma on higher moments of \ $(X_n)_{n\geq -1}$.

\begin{Lem}\label{Lem_seged_momentr}
Let \ $(X_n)_{ n\geq -1}$ \ be a second-order Galton--Watson process (without immigration) such
 that \ $\EE(X_{-1}^r) < \infty$, \ $\EE(X_0^r) < \infty$, \ $\EE(\xi^r) < \infty$ \ and
 \ $\EE(\eta^r) < \infty$ \ with some \ $r > 1$.
\ Then \ $\EE(X_n^r) < \infty$ \ for all \ $n \in \NN$.
\end{Lem}

\noindent{\bf Proof.}
By power means inequality, we have
 \begin{align*}
  \EE(X_n^r \mid \cF_{n-1}^X)
  &= \EE\left(\left(\sum_{i=1}^{X_{n-1}} \xi_{n,i}
                    + \sum_{j=1}^{X_{n-2}} \eta_{n,j}\right)^r \,\Bigg|\, \cF_{n-1}^X\right) \\
  &\leq 2^{r-1}
        \EE\left(\left(\sum_{i=1}^{X_{n-1}} \xi_{n,i}\right)^r
                 + \left(\sum_{j=1}^{X_{n-2}} \eta_{n,j}\right)^r \,\Bigg|\ \cF_{n-1}^X\right)\\
  &\leq 2^{r-1}
        \EE\left(X_{n-1}^{r-1} \sum_{i=1}^{X_{n-1}} \xi_{n,i}^r
                 + X_{n-2}^{r-1} \sum_{j=1}^{X_{n-2}} \eta_{n,j}^r \,\Bigg|\ \cF_{n-1}^X\right) \\
  &= 2^{r-1}
     \left(X_{n-1}^r \EE(\xi^r) + X_{n-2}^r \EE(\eta^r)\right)
   < \infty
 \end{align*}
 for all \ $n \in \NN$.
\ Hence \ $\EE(X_n^r) \leq 2^{r-1} \big(\EE(X_{n-1}^r)\EE(\xi^r) +  \EE(X_{n-2}^r)\EE(\eta^r)\big)$,
 \ $n \in \NN$.
\ By induction we obtain the statement.
\proofend

Moreover, we present an auxiliary lemma on an estimation of the second moment of
 a second-order Galton--Watson process (without immigration).
This lemma is valid for the subcritical, critical and supercritical cases as well,
 however, in the proofs we only use it for the subcritical case.

\begin{Lem}\label{Lem_seged_moment}
Let \ $(X_n)_{ n\geq -1}$ \ be a second-order Galton--Watson process (without immigration) such
 that \ $X_0 = 1$, \ $X_{-1} = 0$, \ $\EE(\xi^2) < \infty$ \ and \ $\EE(\eta^2) < \infty$.
\ Then for all \ $n \in \NN$,
 \begin{align}\label{help8}
  \EE(X_n^2)
  \leq \begin{cases}
        c_\SUB \, \varrho^n , & \text{if \ $\varrho \in (0, 1)$,} \\
        c_\CRIT \, n , & \text{if \ $\varrho = 1$,} \\
        c_\SUP \, \varrho^{2n} , & \text{if \ $\varrho \in (1, \infty)$,}
       \end{cases}
 \end{align}
 where
 \[
   c_\SUB := 1 + \frac{\var(\xi)}{\varrho(1-\varrho)}
             + \frac{\var(\eta)}{\varrho^2(1-\varrho)} , \quad
   c_\CRIT := 1 + \var(\xi) + \var(\eta) , \quad
   c_\SUP := 1 + \frac{\var(\xi)}{\varrho(\varrho-1)}
             + \frac{\var(\eta)}{\varrho^3(\varrho-1)} .
 \]
\end{Lem}

\noindent{\bf Proof.}
By formula (A2) in Lemma A.1 in Isp\'any and Pap \cite{IspPap}, we have
 \[
   \var(\bY_n)
   = \sum_{j=0}^{n-1}
      \bM_{\xi,\eta}^j
      \bigl[(\be_1^\top  \bM_{\xi,\eta}^{n-j-1}  \be_1) \bV_\xi
            + (\be_2^\top \bM_{\xi,\eta}^{n-j-1} \be_1) \bV_\eta\bigr]
      ( \bM_{\xi,\eta}^\top )^j ,
 \]
 where  \ $(\bY_n)_{ n\in\ZZ_+}$ \ is given by \eqref{help15}  with \ $\bY_0 = [1\; 0]^\top$, \ and
 \[
   \bV_\xi
   := \var\left(\begin{bmatrix} \xi \\ 1 \end{bmatrix}\right)
   = \var(\xi) \be_1 \be_1^\top , \qquad
   \bV_\eta := \var\left(\begin{bmatrix} \eta \\ 0 \end{bmatrix}\right)
   = \var(\eta) \be_1 \be_1^\top ,
 \]
 where \ $\xi$ \ and \ $\eta$ \ are random variables such that \ $\xi\distre \xi_{1,1}$ \ and
 \ $\eta\distre \eta_{1,1}$.
\ Here we note that formula (A2) in Lemma A.1 in Isp\'any and Pap \cite{IspPap} is stated  only for
 critical processes, but it also holds in the subcritical and supercritical cases as well; the
 proof is the very same.
Consequently,
 \begin{align*}
  \var(X_n)
  &= \var(\be_1^\top \bY_n) = \be_1^\top \var(\bY_n) \be_1 \\
  &= \be_1^\top
     \sum_{j=0}^{n-1}
      \bM_{\xi,\eta}^j
      \bigl[(\be_1^\top \bM_{\xi,\eta}^{n-j-1} \be_1) \var(\xi) \be_1 \be_1^\top
            + (\be_2^\top \bM_{\xi,\eta}^{n-j-1} \be_1) \var(\eta) \be_1 \be_1^\top\bigr]
      ( \bM_{\xi,\eta}^\top )^j
      \be_1 \\
  &= \sum_{j=0}^{n-1}
      (\be_1^\top \bM_{\xi,\eta}^j \be_1)^2
      \bigl[ \var(\xi) (\be_1^\top \bM_{\xi,\eta}^{n-j-1} \be_1)
            + \var(\eta) (\be_2^\top \bM_{\xi,\eta}^{n-j-1} \be_1)\bigr] ,
 \end{align*}
where we used that
 \ $\be_1^\top (\bM_{\xi,\eta}^\top)^j \be_1 = \be_1^\top \bM_{\xi,\eta}^j \be_1$.
\ Using \eqref{EX} with \ $X_0 = 1$ \ and \ $X_{-1} = 0$, \ we have
 \ $\be_1^\top \bM_{\xi,\eta}^j \be_1 = \EE(X_j)$ \ and
 \ $\be_2^\top \bM_{\xi,\eta}^j \be_1 = \EE(X_{j-1})$ \ for each \ $j \in \ZZ_+$, \ hence
 \[
   \var(X_n)
   = \var(\xi) \sum_{j=0}^{n-1} [\EE(X_j)]^2 \EE(X_{n-j-1})
     + \var(\eta) \sum_{j=0}^{n-2} [\EE(X_j)]^2 \EE(X_{n-j-2}) ,
 \]
 where we used that \ $X_{-1} = 0$.
\ We note that the above formula for \ $\var(X_n)$ \ can also be found in Kashikar and Deshmukh
 \cite[page 562]{KasDes}.
Using \eqref{1moment_est} with \ $X_0=1$ \ and \ $X_{-1}=0$, \ we obtain
 \begin{align*}
  \EE(X_n^2) = \var(X_n) + [\EE(X_n)]^2
  &\leq \var(\xi) \sum_{j=0}^{n-1} \varrho^{n+j-1}
        + \var(\eta) \sum_{j=0}^{n-2} \varrho^{n+j-2} + \varrho^{2n} \\
  &= \begin{cases}
      n \var(\xi) + (n - 1) \var(\eta) + 1 , & \text{if \ $\varrho = 1$,} \\
      \var(\xi) \, \frac{\varrho^{n-1}-\varrho^{2n-1}}{1-\varrho}
      + \var(\eta) \, \frac{\varrho^{n-2}-\varrho^{2n-3}}{1-\varrho} + \varrho^{2n} ,
       & \text{if \ $\varrho \ne 1$,}
     \end{cases}
 \end{align*}
 yielding \eqref{help8}.
Indeed, for example, if \ $\varrho \in (1,\infty)$, \ then
 \begin{align*}
  \frac{\varrho^{n-2} - \varrho^{2n-3}}{1-\varrho}
  = \frac{\varrho^{2n-3} (1-\varrho^{-n+1})}{\varrho -1}
  \leq \frac{\varrho^{2n-3}}{\varrho-1}
  = \frac{\varrho^{2n}}{\varrho^3(\varrho-1)} ,
  \qquad n \in \NN .
 \end{align*}
\proofend

\section{Representation of the unique stationary distribution for 2-type Galton--Watson
         processes with immigration}
\label{App2multitype}

First, we introduce 2-type Galton--Watson processes with immigration.
For each \ $k, j \in \ZZ_+$ \ and \ $i, \ell \in \{ 1, 2 \}$, \ the number of individuals of type
 \ $i$ \ born or arrived as immigrants in the \ $k^\mathrm{th}$ \ generation will be denoted by
 \ $X_{k,i}$, \ the number of type \ $\ell$ \ offsprings produced by the \ $j^\mathrm{th}$
 \ individual who is of type \ $i$ \ belonging to the \ $(k-1)^\mathrm{th}$ \ generation
 will be denoted by \ $\xi_{k,j,i,\ell}$, \ and the number of type \ $i$ \ immigrants in
 the \ $k^\mathrm{th}$ \ generation will be denoted by \ $\vare_{k,i}$.
\ Then we have
 \begin{equation}\label{GWI(2)}
  \begin{bmatrix}
   X_{k,1} \\
   X_{k,2}
  \end{bmatrix}
  = \sum_{j=1}^{X_{k-1,1}}
     \begin{bmatrix}
      \xi_{k,j,1,1} \\
      \xi_{k,j,1,2}
     \end{bmatrix}
    +\sum_{j=1}^{X_{k-1,2}}
       \begin{bmatrix}
        \xi_{k,j,2,1} \\
        \xi_{k,j,2,2}
       \end{bmatrix}
    + \begin{bmatrix}
       \vare_{k,1} \\
       \vare_{k,2}
      \end{bmatrix} ,
  \qquad k \in \NN .
 \end{equation}
Here
 \ $\bigl\{\bX_0, \, \bxi_{k,j,i}, \, \bvare_k : k, j \in \NN, \, i \in \{1, 2\}\bigr\}$
 \ are supposed to be independent, and \ $\{\bxi_{k,j,1} : k, j \in \NN\}$,
 \ $\{\bxi_{k,j,2} : k, j \in \NN\}$ \ and \ $\{\bvare_k : k \in \NN\}$ \ are supposed to
 consist of identically distributed random vectors, where
 \[
   \bX_0 := \begin{bmatrix}
             X_{0,1} \\
             X_{0,2}
            \end{bmatrix} , \qquad
   \bX_k := \begin{bmatrix}
             X_{k,1} \\
             X_{k,2}
            \end{bmatrix} , \qquad
   \bxi_{k,j,i} := \begin{bmatrix}
                  \xi_{k,j,i,1} \\
                  \xi_{k,j,i,2}
                 \end{bmatrix} , \qquad
   \bvare_k := \begin{bmatrix}
                \vare_{k,1} \\
                \vare_{k,2}
               \end{bmatrix} .
 \]
For notational convenience, let \ $\bxi_1$, \ $\bxi_2$ \ and \ $\bvare$ \ be random vectors such
 that \ $\bxi_1 \distre \bxi_{1,1,1}$, \ $\bxi_2 \distre \bxi_{1,1,2}$ \ and
 \ $\bvare \distre \bvare_1$, \ and put \ $\bm_{\bxi_1} := \EE(\bxi_1) \in [0, \infty]^2$,
 \ $\bm_{\bxi_2} := \EE(\bxi_2) \in [0, \infty]^2$, \ and
 \ $\bm_\bvare := \EE(\bvare) \in [0, \infty]^2$, \ and put
 \[
   \bM_{\bxi} := \begin{bmatrix}
                  \bm_{\bxi_1} & \bm_{\bxi_2}
                 \end{bmatrix} \in [0, \infty]^{2 \times 2} .
 \]
We call \ $\bM_{\bxi}$ \ the offspring mean matrix, and note that many authors define the offspring
 mean matrix as \ $\bM^\top_\bxi$.
\ If \ $\bm_{\bxi_1} \in \RR^2_+$, \ $\bm_{\bxi_2} \in \RR^2_+$, \ and \ $\bm_\bvare \in \RR^2_+$,
 \ then for each \ $n \in \ZZ_+$, \ \eqref{GWI(2)} implies
 \begin{equation*}
  \EE(\bX_n \mid \cF_{n-1}^\bX)
  = X_{n-1,1} \, \bm_{\bxi_1} + X_{n-1,2} \, \bm_{\bxi_2} + \bm_\bvare
  = \bM_{\bxi} \, \bX_{n-1} + \bm_\bvare , \qquad n \in \NN ,
 \end{equation*}
 where \ $\cF_n^\bX := \sigma\bigl(\bX_0, \dots, \bX_n\bigr)$, \ $n \in \ZZ_+$.
\ Consequently, \ $\EE(\bX_n) = \bM_{\bxi} \EE(\bX_{n-1}) + \bm_\bvare$, \ $n \in \NN$, \ which
 implies
 \begin{equation*}
  \EE(\bX_n) = \bM_{\bxi}^n \, \EE(\bX_0) + \sum_{k=1}^n \bM_{\bxi}^{n-k} \bm_\bvare , \qquad
  n \in \NN .
 \end{equation*}
Hence, the asymptotic behavior of the sequence \ $(\EE(\bX_n))_{n\in\ZZ_+}$ \ depends on the
 asymptotic behavior of the powers \ $(\bM_{\bxi}^n)_{n\in\NN}$ \ of the offspring mean matrix,
 which is related to the spectral radius \ $r(\bM_\bxi) \in \RR_+$ \ of \ $\bM_\bxi$ \ (see the
 Frobenius--Perron theorem, e.g., Horn and Johnson \cite[Theorems 8.2.8 and 8.5.1]{HorJoh}).
A 2-type Galton--Watson process \ $(\bX_n)_{n\in\ZZ_+}$ \ with immigration is referred to
 respectively as \emph{subcritical}, \emph{critical} or \emph{supercritical} if
 \ $r(\bM_\bxi) < 1$, \ $r(\bM_\bxi) = 1$ \ or \ $r(\bM_\bxi) > 1$ \ (see, e.g., Athreya and Ney
 \cite[V.3]{AthNey} or Quine \cite{Qui}).
We extend this classification for all 2-type Galton--Watson processes with immigration.

If \ $\bm_{\bxi_1} \in \RR_+^2$, \ $\bm_{\bxi_2} \in \RR_+^2$, \ $r(\bM_\bxi) < 1$,
 \ $\bM_\bxi$ \ is primitive, i.e., there exists \ $m \in \NN$ \ such that
 \ $\bM_\bxi^m \in \RR^{2 \times 2}_{++}$, \ $\PP(\bvare = \bzero) < 1$ \ and
 \ $\EE(\bbone_{\{\bvare\ne\bzero\}} \log((\be_1+\be_2)^\top\bvare)) < \infty$, \ then, by the
 Theorem in Quine \cite{Qui}, there exists a unique stationary distribution \ $\bpi$ \ for
 \ $(\bX_n)_{n\in\ZZ_+}$.
\ As a consequence of formula (16) for the probability generating function of \ $\bpi$ \ in Quine
 \cite{Qui}, we have
 \[
   \sum_{i=0}^n \bV_i^{(i)}(\bvare_i) \distr \bpi \qquad \text{as \ $n \to \infty$,}
 \]
 where \ $(\bV_k^{(i)}(\bvare_i))_{k\in\ZZ_+}$, \ $i \in \ZZ_+$, \ are independent copies of a
 2-type Galton--Watson process \ $(\bV_k(\bvare))_{k\in\ZZ_+}$ \ (without immigration) with initial
 vector \ $\bV_0(\bvare) = \bvare$ \ and with the same offspring distributions as
 \ $(\bX_k)_{k\in\ZZ_+}$.
\ Consequently, we have
 \begin{equation*}
  \sum_{i=0}^\infty \bV_i^{(i)}(\bvare_i) \distre \bpi ,
 \end{equation*}
 where the series \ $\sum_{i=0}^\infty \bV_i^{(i)}(\bvare_i)$ \ converges with probability 1, see,
 e.g., Heyer \cite[Theorem 3.1.6]{Hey}.
The above representation of the stationary distribution \ $\bpi$ \ for \ $(\bX_n)_{n\in\ZZ_+}$ \
 can be interpreted in a way that we consider independent 2-type Galton--Watson processes without immigration such that
 the \ $i^\mathrm{th}$ \ one admits initial vector \ $\bvare_i$, \ $i\in\ZZ_+$, \ evaluate the \ $i^\mathrm{th}$ \ 2-type Galton-Watson processes at time point \ $i$, \
 and then sum up all these random variables.

\section{Regularly varying distributions}
\label{App1}

First, we recall the notions of slowly varying and regularly varying functions, respectively.

\begin{Def}
A measurable function \ $U: \RR_{++} \to \RR_{++}$ \ is called regularly varying at infinity with
 index \ $\rho \in \RR$ \ if for all \ $q \in \RR_{++}$,
  \[
    \lim_{x\to\infty} \frac{U(qx)}{U(x)} = q^\rho .
  \]
In case of \ $\rho = 0$, \ $U$ \ is called slowly varying at infinity.
\end{Def}

Next, we recall the notion of regularly varying random variables.

\begin{Def}
A non-negative random variable \ $X$ \ is called regularly varying with index \ $\alpha \in \RR_+$
 \ if \ $U(x) := \PP(X > x) \in \RR_{++}$ \ for all \ $x \in \RR_{++}$, \ and \ $U$ \ is regularly
 varying at infinity with index \ $-\alpha$.
\end{Def}

\begin{Lem}\label{Lem_regl_power}
If \ $\zeta$ \ is a non-negative regularly varying random variable with index \ $\alpha \in \RR_+$,
 \ then for each \ $c \in \RR_{++}$, \ $\zeta^c$ \ is regularly varying with index
 \ $\frac{\alpha}{c}$.
\end{Lem}

\noindent{\bf Proof.}
For any \ $q \in \RR_{++}$, \ we have
 \[
   \lim_{x\to\infty} \frac{\PP(\zeta^c>qx)}{\PP(\zeta^c>x)}
   = \lim_{x\to\infty} \frac{\PP(\zeta>q^{1/c}x^{1/c})}{\PP(\zeta>x^{1/c})}
   = q^{-\alpha/c} ,
 \]
 as desired.
\proofend

\begin{Lem}\label{sv}
If \ $L : \RR_{++} \to \RR_{++}$ \ is a slowly varying function (at infinity), then
 \[
   \lim_{x\to\infty} x^\delta L(x) = \infty , \qquad
   \lim_{x\to\infty} x^{-\delta} L(x) = 0 , \qquad \delta \in \RR_{++} .
 \]
\end{Lem}
For Lemma \ref{sv}, see, Bingham et al.\ \cite[Proposition 1.3.6. (v)]{BinGolTeu}.

\begin{Lem}\label{help_log_exp}
If \ $\vare$ \ is a non-negative regularly varying random variable with index
 \ $\alpha \in \RR_{++}$, \ then \ $\EE(\bbone_{\{\vare \ne 0\}} \log(\vare)) < \infty$ \ and
 \ $\EE(\log(\vare + 1)) < \infty$.
\end{Lem}

\noindent{\bf Proof.}
Since \ $\EE(\bbone_{\{\vare \ne 0\}} \log(\vare)) \leq \EE(\log(\vare + 1))$, \ it is enough to prove
 that \ $\EE(\log(\vare + 1)) < \infty$.
\ Since \ $\log(\vare + 1)\geq 0$, \ we have
 \begin{align*}
  \EE(\log(\vare + 1))
  &= \int_0^\infty \PP(\log(\vare + 1) \geq x) \, \dd x
   = \int_0^\infty \PP(\vare \geq \ee^x - 1) \, \dd x \\
  &= \int_0^1 \PP(\vare \geq \ee^x - 1) \, \dd x
     + \int_1^\infty \PP(\vare \geq \ee^x - 1) \, \dd x
   := I_1 + I_2 .
 \end{align*}
Here \ $I_1 \leq 1$, \ and, by substitution \ $y = \ee^x - 1$,
 \[
   I_2 = \int_{\ee-1}^\infty y^{-\alpha} L(y) \frac{1}{1+y} \, \dd y ,
 \]
 where \ $L(y) := y^\alpha \PP(\vare > y)$, \ $y \in \RR_{++}$, \ is a slowly varying function.
By Lemma \ref{sv}, there exists \ $y_0 \in (\ee - 1, \infty)$ \ such that
 \ $y^{-\frac{\alpha}{2}} L(y) \leq 1$ \ for all \ $y \in [y_0, \infty)$.
\ Hence
 \begin{align*}
  I_2 &= \int_{\ee-1}^{y_0} y^{-\alpha} L(y) \frac{1}{1+y} \, \dd y
         + \int_{y_0}^\infty y^{-\alpha} L(y) \frac{1}{1+y} \, \dd y \\
      &\leq \int_{\ee-1}^{y_0} y^{-\alpha} L(y) \frac{1}{1+y} \, \dd y
            + \int_{y_0}^\infty y^{-\frac{\alpha}{2}} \frac{1}{1+y} \, \dd y \\
      &\leq \int_{\ee-1}^{y_0} y^{-\alpha} L(y) \frac{1}{1+y} \, \dd y
            + \int_{y_0}^\infty y^{-\frac{\alpha}{2}-1} \, \dd y \\
      &\leq \int_{\ee-1}^{y_0} \frac{1}{1+y} \, \dd y
            + \int_{y_0}^\infty y^{-\frac{\alpha}{2}-1} \, \dd y
       < \infty ,
 \end{align*}
 since \ $y^{-\alpha} L(y) = \PP(\vare > y) \leq 1$ \ for all \ $y \in \RR_{++}$.
\proofend

\begin{Lem}\label{Lem_regl_inequality}
If \ $\eta$ \ is a non-negative regularly varying random variable with index
 \ $\alpha \in(1, 2)$, \ then for every \ $\varrho \in (\alpha, \infty)$, \ there exist
 \ $y_0 \in \RR_{++}$ \ and \ $B \in \RR_{++}$ \ such that
 \[
   \frac{\PP(\eta>z)}{\PP(\eta>y)} \leq B \left(\frac{z}{y}\right)^{-\varrho} ,
   \qquad y \geq z \geq y_0 ,
 \]
 or equivalently,
 \[
   \frac{\PP(\eta>\theta y)}{\PP(\eta>y)} \leq B \theta^{-\varrho} ,
   \qquad \theta\in(0,1], \qquad y\geq \frac{y_0}{\theta} .
 \]
\end{Lem}
For Lemma \ref{Lem_regl_inequality}, see Proposition 2.2.1 in Bingham et al.\ \cite{BinGolTeu}.

For the next lemma, see Fa\"{y} et al.\ \cite[Lemma 4.4]{GilGAMikSam}.
Here we present a proof as well, since we state their result in a little bit extended form.

\begin{Lem}\label{Lem_L_const}
Let \ $h : \RR_+ \to \RR_{++}$ \ be a function such that \ $\lim_{x\to\infty} h(x) = 0$.
\ Then there exists a monotone increasing, left-continuous, slowly varying (at infinity) function
 \ $L$ \ such that \ $L(x) \geq 1$, \ $x \in \RR_+$, \ $\lim_{x\to\infty} L(x) = \infty$ \ and
 \ $\lim_{x\to\infty} L(x) h(x) = 0$.
\ One can also choose a version of \ $L$ \ which is right-continuous with all the other properties
 remaining true.
\end{Lem}

\noindent{\bf Proof.}
We can construct \ $L$ \ as follows.
Let \ $L(x) := 1$ \ for \ $x \in [0, x_0]$, \ where \ $x_0 := \sup\{y \in \RR_+ : h(y) >1\}$, \ and
 we define \ $\sup\emptyset := 0$.
\ Since \ $\lim_{x\to\infty} h(x) = 0$, \ we have \ $x_0 \in \RR_+$.
\ Let \ $L(x) := 2$ \ for \ $x \in (x_0, x_1]$, \ where
 \ $x_1 := \max\{2 x_0, \sup\{y \in \RR_+ : h(y) > 2^{-2}\}\}$.
\ Let \ $L(x) := 3$ \ for \ $x \in (x_1, x_2]$, \ where
 \ $x_2 := \max\{3 x_1, \sup\{y \in \RR_+ : h(y) > 3^{-2}\}\}$, \ and continue this construction in
 the straightforward way: \ $L(x) := k + 1$ \ for \ $x \in (x_{k-1}, x_k]$, \ where
 \ $x_k := \max\{(k + 1) x_{k-1}, \sup\{y \in \RR_+ : h(y) > (k + 1)^{-2}\}\}$, \ $k \in \NN$.
\ Since \ $h$ \ takes positive values and \ $\lim_{x\to\infty} h(x) = 0$, \ we have
 \ $\lim_{x\to\infty} L(x) = \infty$, \ and, since for all \ $k \in \ZZ_+$ \ and \ $x > x_k$,
 \[
   L(x) h(x)
   = \sum_{i=k}^\infty L(x)h(x)\bbone_{(x_i,x_{i+1}](x)}
   \leq \sum_{i=k}^\infty (i + 2) \frac{1}{(i+1)^2} \bbone_{(x_i,x_{i+1}](x)}
   \leq \frac{k+2}{(k+1)^2} ,
 \]
 we have \ $\lim_{x\to\infty} L(x) h(x) = 0$.
\ It remains to check that \ $L$ \ is slowly varying (at infinity).
For this it is enough to verify that for any \ $ q \in \RR_{++}$ \ and sufficiently large
 \ $ x\in\RR_{++}$, \ we have \ $x$ \ and \ $qx$ \ are either in the same interval of type
 \ $(x_{k-1}, x_k]$ \ or in two neighbouring intervals of this type, since in this case for sufficiently large \ $ x \in \RR_{++}$:
 \[
   \frac{L(qx)}{L(x)} \in \left\{1, \frac{k_x}{k_x+1}, \frac{k_x+1}{k_x}\right\}
 \]
 with some \ $k_x \in \NN$, \ and for sufficiently large \ $ x \in \RR_{++}$ \ and for
 \ $y \geq x$,
 \begin{align*}
  \left|\frac{L(qy)}{L(y)} - 1\right|
  \in \left\{0, \left|\frac{k_y+1}{k_y} - 1\right|, \left|\frac{k_y}{k_y+1} - 1\right|\right\}
      = \left\{0, \frac{1}{k_y}, \frac{1}{k_y+1}\right\} ,
 \end{align*}
 where \ $k_y \geq k_x$ \ and \ $\lim_{y\to\infty} k_y = \infty$, \ yielding that
 \ $\lim_{x\to\infty} \frac{L(qx)}{L(x)} = 1$.
\ To finish the proof, if \ $x \in (x_{k-1}, x_k]$ \ with some \ $k \in \NN$, \ then in case of
 \ $q \geq 1$, \ we have \ $q x \in (x_{k-1}, x_k] \cup (x_k, x_{k+1}]$ \ provided that
 \ $k + 2 \geq q$, \ and in case of \ $q \in (0, 1)$, \ we have
 \ $q x \in (x_{k-2}, x_{k-1}] \cup (x_{k-1}, x_k]$ \ provided that \ $k > \frac{1}{q}$.
\ Indeed, \ $x_k \geq (k + 1) x_{k-1}$, \ $k \in \NN$, \ and if \ $k + 2 \geq q$, \ then
 \ $q x \leq q x_k \leq (k + 2) x_k \leq x_{k+1}$, \ as desired, and if \ $k > \frac{1}{q}$, \ then
 \ $q x > q x_{k-1} > \frac{1}{k} x_{k-1} \geq x_{k-2}$, \ as desired.
\proofend

\begin{Lem}\label{exposv}
If \ $X$ \ and \ $Y$ \ are non-negative random variables such that \ $X$ \ is regularly varying
 with index \ $\alpha \in \RR_+$ \ and there exists \ $r \in (\alpha, \infty)$ \ with
 \ $\EE(Y^r) < \infty$, \ then \ $\PP(Y > x) = \oo(\PP(X > x))$ \ as \ $x \to \infty$.
\end{Lem}
For a proof of Lemma \ref{exposv}, see, e.g., Barczy et al.\ \cite[Lemma C.6]{BarBosPap}.

\begin{Lem}\label{svsv}
If \ $X_1$ \ and \ $X_2$ \ are non-negative regularly varying random variables with index
 \ $\alpha_1 \in \RR_+$ \ and \ $\alpha_2 \in \RR_+$, \ respectively, such that
 \ $\alpha_1 < \alpha_2$, \ then \ $\PP(X_2 > x) = \oo(\PP(X_1 > x))$ \ as \ $x \to \infty$.
\end{Lem}
For a proof of Lemma \ref{svsv}, see, e.g., Barczy et al.\ \cite[Lemma C.7]{BarBosPap}.

\begin{Lem}[Convolution property]\label{Lem_conv}
If \ $X_1$ \ and \ $X_2$ \ are non-negative random variables such that \ $X_1$ \ is regularly varying
 with index \ $\alpha_1 \in \RR_+$ \ and \ $\PP(X_2 > x) = \oo(\PP(X_1 > x))$ \ as
 \ $x \to \infty$, \ then \ $\PP(X_1 + X_2 > x) \sim \PP(X_1 > x)$ \ as \ $x \to \infty$, \ and
 hence \ $X_1 + X_2$ \ is regularly varying with index \ $\alpha_1$.

If \ $X_1$ \ and \ $X_2$ \ are independent non-negative regularly varying random variables with
 index \ $\alpha_1 \in \RR_+$ \ and \ $\alpha_2 \in \RR_+$, \ respectively, then
 \[
   \PP(X_1 + X_2 > x)
   \sim \begin{cases}
         \PP(X_1 > x) & \text{if \ $\alpha_1 < \alpha_2$,} \\
         \PP(X_1 > x) + \PP(X_2 > x) & \text{if \ $\alpha_1 = \alpha_2$,} \\
         \PP(X_2 > x) & \text{if \ $\alpha_1 > \alpha_2$,}
        \end{cases}
 \]
 as \ $x \to \infty$, \ and hence \ $X_1 + X_2$ \ is regularly varying with index
 \ $\min\{\alpha_1, \alpha_2\}$.
\end{Lem}
The statements of Lemma \ref{Lem_conv} follow, e.g., from parts 1 and 3 of Lemma B.6.1 of
 Buraczewski et al.\ \cite{BurDamMik} and Lemma \ref{svsv} together with the fact that the sum of
 two slowly varying functions is slowly varying.

\begin{Thm}[Karamata's theorem]\label{Krthm}
Let \ $U : \RR_{++} \to \RR_{++}$ \ be a locally integrable function such that it is integrable on
 intervals including \ $0$ \ as well.

\noindent
(i) If \ $U$ \ is regularly varying (at infinity) with index \ $-\alpha \in [-1, \infty)$, \ then
 \ $ \RR_{++} \ni x \mapsto \int_0^x U(t) \, \dd t$ \ is regularly varying (at infinity) with index
 \ $1-\alpha$, \ and
 \[
   \lim_{x\to\infty} \frac{xU(x)}{\int_0^x U(t)\,\dd t} = 1 - \alpha .
 \]
\noindent
(ii) If \ $U$ \ is regularly varying (at infinity) with index \ $-\alpha \in (-\infty,-1)$, \ then
 \ $\RR_{++} \ni x \mapsto \int_x^\infty U(t) \, \dd t$ \ is regularly varying (at infinity) with
 index \ $1-\alpha$, \ and
 \[
   \lim_{x\to\infty} \frac{xU(x)}{\int_x^\infty U(t)\,\dd t} = - 1 +\alpha .
 \]
\end{Thm}
For Theorem \ref{Krthm}, see, e.g., Resnick \cite[Theorem 2.1]{Res}.

\begin{Lem}[Potter's bounds]\label{Pb}
If \ $U :  \RR_{++} \to \RR_{++}$ \ is a regularly varying function (at infinity)
 with index \ $  -\alpha \in \RR$, \ then for every \ $\delta  \in\RR_{++}$, \ there
 exists \ $x_0 \in \RR_+$ \ such that
 \[
   (1 - \delta) q^{ -\alpha - \delta}
   < \frac{U(q x)}{U(x)}
   < (1 + \delta) q^{-\alpha + \delta} ,
   \qquad x \in [x_0, \infty) , \quad q \in [1, \infty) .
 \]
\end{Lem}
For Lemma \ref{Pb}, see, e.g., Resnick \cite[Proposition 2.6]{Res}.

Finally, we recall a result on the tail behaviour of regularly varying random sums.

\begin{Pro}\label{FGAMSRS}
Let \ $\tau$ \ be a non-negative integer-valued random variable and let
 \ $\{\zeta, \zeta_i : i \in \NN\}$ \ be independent and identically distributed non-negative
 random variables, independent of \ $\tau$, \ such that \ $\tau$ \ is regularly varying with index
 \ $\beta \in \RR_+$ \ and \ $\EE(\zeta) \in \RR_{++}$.
\ In case of \ $\beta \in [1,\infty)$, \ assume additionally that there exists
 \ $r \in (\beta, \infty)$ \ with \ $\EE(\zeta^r) < \infty$.
\ Then we have
 \[
   \PP\biggl(\sum_{i=1}^\tau \zeta_i > x\biggr)
   \sim \PP\biggl(\tau > \frac{x}{\EE(\zeta)}\biggr)
   \sim (\EE(\zeta))^\beta \PP(\tau > x) \qquad \text{as \ $x \to \infty$,}
 \]
 and hence \ $\sum_{i=1}^\tau \zeta_i$ \ is also regularly varying with index \ $\beta$.
\end{Pro}

For a proof of Proposition \ref{FGAMSRS}, see, e.g.,  Barczy et al.\
 \cite[Proposition F.3]{BarBosPap}.

\section{Large deviations}
\label{App3}

We recall a result about large deviations for sums of non-negative independent and identically
 distributed regularly varying random variables, see, Tang and Yan
 \cite[part (ii) of Theorem 1]{TangYan}.
We use it in the second proof of Theorem \ref{2GW_X_0_X_-1} in case of $\alpha\in(1,2)$.
Here we present a complete proof as well, since the one in Tang and Yan \cite[part (ii) of Theorem 1]{TangYan} contains a gap.

\begin{Thm}[Large deviations]\label{Thm_largedev}
If \ $(\eta_j)_{j\in\NN}$ \ are independent, identically distributed non-negative regularly
 varying random variables with index \ $\alpha \in (1, 2)$, \ then for each
 \ $\gamma \in (\EE(\eta_1), \infty)$, \ there exists a constant \ $C \in \RR_{++}$ \ such that
 \[
   \PP(\eta_1 + \cdots + \eta_n > y) \leq C n \PP(\eta_1 > y)
 \]
 for all \ $n \in \NN$ \ and \ $y \in [\gamma n, \infty)$.
\end{Thm}

\noindent{\bf Proof.}
We will follow the proof of part (ii) of Theorem 1 in Tang and Yan \cite{TangYan}.
Let \ $q \in (0, 1)$ \ and
 \[
   \widetilde\eta_j := \eta_j\bbone_{\{\eta_j\leq qy\}}, \quad j \in \NN ,
   \qquad \tS_j := \sum_{i=1}^j \widetilde\eta_i , \quad j \in \NN.
 \]
Then for all \ $n \in \NN$,
 \begin{align}\label{help9_large}
  \begin{split}
   &\PP(\eta_1 + \cdots + \eta_n > y) \\
   &\qquad = \PP(\eta_1 + \cdots + \eta_n > y, \max_{1\leq j\leq n} \eta_j > q y)
             + \PP(\eta_1 + \cdots + \eta_n > y, \max_{1\leq j\leq n} \eta_j \leq q y) \\
   &\qquad \leq \PP(\max_{1\leq j\leq n} \eta_j > q y)
                + \PP(\widetilde\eta_1 + \cdots + \widetilde\eta_n > y,
                      \max_{1\leq j\leq n} \eta_j \leq q y) \\
   &\qquad \leq \PP(\max_{1\leq j\leq n} \eta_j > q y) + \PP(\tS_n > y) \\
   &\qquad \leq \sum_{j=1}^n \PP(\eta_j > q y) + \PP(\tS_n > y) \\
   &\qquad = n\PP(\eta_1 > q y) + \PP(\tS_n > y) ,
    \qquad y \in \RR_+ .
  \end{split}
 \end{align}
Here
 \[
   \PP(\eta_1 > q y)
   = \frac{\PP(\eta_1>qy)}{\PP(\eta_1>y)} \cdot \PP(\eta_1 > y) , \qquad y \in \RR_{++} ,
 \]
 and since \ $\lim_{y\to\infty} \frac{\PP(\eta_1>qy)}{\PP(\eta_1>y)} = q^{-\alpha}$, \ there exists
 an \ $y_* \in \RR_{++}$ \ such that \ $\frac{\PP(\eta_1>qy)}{\PP(\eta_1>y)} \leq 2 q^{-\alpha}$
 \ for all \ $y \geq y_*$.
\ Now we check that \ $\frac{\PP(\eta_1>qy)}{\PP(\eta_1>y)}$ \ is bounded on the interval
 \ $[0, y_*]$.
\ Since \ $\lim_{y\to 0} \frac{\PP(\eta_1>qy)}{\PP(\eta_1>y)} = 1$, \ there exists an
 \ $y_1 \in \RR_{++}$ \ such that \ $y_1 < y_*$ \ and
 \ $\frac{\PP(\eta_1>qy)}{\PP(\eta_1>y)} \leq 2$ \ on the interval \ $[0, y_1]$.
\ On the interval \ $[y_1, y_{*}]$ \ the quantity \ $\frac{\PP(\eta_1>qy)}{\PP(\eta_1>y)}$
 \ can be bounded from above by \ $\frac{\PP(\eta_1>qy_1)}{\PP(\eta_1>y_{*})}$.
\ Hence the function \ $\RR_+ \ni y \mapsto \frac{\PP(\eta_1>qy)}{\PP(\eta_1>y)}$ \ is bounded, and
 consequently, there exists a constant \ $C_1(q) \in \RR_{++}$ \ (depending possibly on the
 distribution of \ $\eta_1$ \ as well) such that
 \begin{align}\label{help10_large}
  n \PP(\eta_1 > q y) \leq C_1(q) n \PP(\eta_1 > y) , \qquad y \in \RR_{++} , \quad n \in \NN .
 \end{align}
Let \ $a(n, y) := \max\{-\log(n \PP(\eta_1 > y), 1\}$, \ $n \in \NN$, \ $y \in \RR_{++}$.
\ Then \  $a(n, y)$ \ tends to \ $\infty$ \ uniformly for \ $y \geq \gamma n$ \ as
 \ $n \to \infty$, \ i.e., \ $\lim_{n\to\infty} \inf_{y\geq\gamma n} a(n, y) = \infty$, \ since, by
 Lemma \ref{sv},
 \begin{align}\label{help20}
  n \PP(\eta_1 > y)
  \leq n \PP(\eta_1 > \gamma n)
  = n (\gamma n)^{-\alpha} L_{\eta_1}(\gamma n)
  = \gamma^{-\alpha} n^{1-\alpha} L_{\eta_1}(n) \frac{L_{\eta_1}(\gamma n)}{L_{\eta_1}(n)}
  \to \gamma^{-\alpha} \cdot 0 \cdot 1 = 0
 \end{align}
 as \ $n \to \infty$, \ where \ $L_{\eta_1}(y) := y^\alpha \PP(\eta_1 > y)$,
 \ $y \in \RR_{++}$, \ is a slowly varying (at infinity) function.
For any \ $y \in \RR_{++}$, \ $h \in \RR_{++}$ \ and \ $n \in \NN$, \ we have
 \begin{align*}
  \frac{\PP(\tS_n>y)}{n\PP(\eta_1>y)}
  &\leq \frac{\ee^{-hy}\EE(\ee^{h\tS_n})}{n\PP(\eta_1>y)}
   = \frac{\ee^{-hy}(\EE(\ee^{h\widetilde\eta_1}))^n}{n\PP(\eta_1>y)} \\
  &= \frac{\ee^{-hy}\left(\int_0^{qy}\ee^{ht}\,\ F_{\eta_1}(\dd t)\right)^n}{n\PP(\eta_1>y)}
   = \frac{\ee^{-hy}\left(\int_0^{qy}(\ee^{ht}-1)\,\ F_{\eta_1}(\dd t)+1\right)^n}
          {n\PP(\eta_1>y)} \\
  &\leq \frac{\ee^{-hy}\exp\left\{n\int_0^{qy}(\ee^{ht}-1)\,F_{\eta_1}(\dd t)\right\}}
             {\ee^{-a(n,y)}} ,
 \end{align*}
 where the last step follows from \ $(1 + y)^n \leq \ee^{ny}$, \ $y \in \RR_+$, \ $n \in \NN$,
 \ and from \ $a(n, y) \geq -\log(n \PP(\eta_1 > y))$, \ yielding
 \ $\ee^{-a(n,y)} \leq n \PP(\eta_1 > y)$.
\ Using that \ $a(n, y) \geq 1$, \ $n \in \NN$, \ $y \in \RR_{++}$, \ let us consider the
 decomposition
 \[
   \int_0^{qy} (\ee^{ht}-1) \, F_{\eta_1}(\dd t)
   = \int_0^{\frac{qy}{a(n,y)}} (\ee^{ht} - 1) \, F_{\eta_1}(\dd t)
     + \int_{\frac{qy}{a(n,y)}}^{qy} (\ee^{ht} - 1) \, F_{\eta_1}(\dd t)
   =: I_1 + I_2 .
 \]
Using the inequality \ $\ee^{y} - 1 \leq y \ee^y$, \ $y \in \RR_+$, \ we have
 \begin{align*}
  I_1 = \int_0^\frac{qy}{a(n,y)} (\ee^{ht} - 1) \, F_{\eta_1}(\dd t)
  &\leq \int_0^{\frac{qy}{a(n,y)}} h t \ee^{ht} \, F_{\eta_1}(\dd t) \\
  &\leq \ee^{\frac{hqy}{a(n,y)}} \int_0^{\frac{qy}{a(n,y)}} h t \, F_{\eta_1}(\dd t)
   \leq h \ee^{\frac{hqy}{a(n,y)}} \EE(\eta_1) .
 \end{align*}
Now we turn to treat \ $I_2$.
\ Applying Lemma \ref{Lem_regl_inequality}, for all \ $\varrho > \alpha$, \ there exist
 \ $y_0 \in \RR_{++}$ \ and \ $B \in \RR_{++}$ \ (possibly depending on \ $\varrho$ \ and on the
 distribution of \ $\eta_1$) \ such that
 \[
   \frac{\PP\left(\eta_1 > \frac{qy}{a(n,y)} \right) }{\PP\left(\eta_1 > y \right)}
   \leq B \left(\frac{q}{a(n,y)}\right)^{-\varrho} \qquad
   \text{whenever \ $y\geq \frac{qy}{a(n,y)} \geq y_0$.}
 \]
The aim of the following discussion is to show that for each \ $n \in \NN$, \ there exists
 \ $\ty_0(n) \in \RR_{++}$ \ such that \ $y \geq \frac{qy}{a(n,y)} \geq y_0$ \ holds for all
 \ $y \geq \ty_0(n)$.
\ For each \ $n \in \NN$, \ the first inequality holds for sufficiently large \ $y$, \ since
 \ $\lim_{y\to\infty} a(n,y) = \infty$.
\ Moreover, for each \ $n \in \NN$, \ the second inequality holds for sufficiently large \ $y$,
 \ since \ $\lim_{y\to\infty} \frac{a(n,y)}{y} = 0$.
\ Indeed, for each \ $n \in \NN$ \ we have \ $a(n,y) = - \log(n \PP(\eta_1 > y))$ \ for
 sufficiently large \ $y$, \ and hence
 \begin{align*}
  \frac{a(n,y)}{y}
  = \frac{-\log(n\PP(\eta_1>y))}{y}
  = \frac{-\log(ny^{-\alpha}L_{\eta_1}(y))}{y}
  = \frac{-\log(n)+\alpha\log(y)-\log(L_{\eta_1}(y))}{y} .
 \end{align*}
By Lemma \ref{sv}, for any \ $\delta \in \RR_{++}$, \ we have
 \ $y^{-\delta} \leq L_{\eta_1}(y) \leq y^\delta$ \ for sufficiently large \ $y$.
\ Taking logarithm, dividing by \ $y$, \ and using that
 \ $\lim_{y\to\infty} \frac{\log(y)}{y} = 0$, \ one concludes
 \ $\lim_{y\to\infty} \frac{a(n,y)}{y} = 0$.

Set
 \[
   h := h(n, y, K) := \frac{a(n,y)-K\varrho\log(a(n,y))}{Kqy} ,
 \]
 where \ $\varrho > \alpha$ \ and \ $K > 1$ \ will be chosen later.
We show that there exists \ $N_1 \in \NN$ \ such that \ $h > 0$ \ and \ $a(n,y) > 1$ \ for all
 \ $y \geq \gamma n$ \ and \ $n \geq N_1$.
\ Since \ $\lim_{x\to\infty} \frac{\log(x)}{x} = 0$, \ there exists \ $M > 0$ \ such that
 \ $\frac{\log(x)}{x} < \frac{1}{K\varrho}$ \ for all \ $x \geq M$.
\ Since \ $\lim_{n\to\infty} \inf_{y\geq \gamma n} a(n,y) = \infty$, \ there exists
 \ $n_0(M) \in \NN$ \ such that \ $a(n,y) \geq M$ \ for all \ $y \geq \gamma n$ \ with
 \ $n \geq n_0(M)$.
\ Hence \ $\frac{\log(a(n,y))}{a(n,y)} < \frac{1}{K\varrho}$ \ for all \ $y \geq \gamma n$ \ with
 \ $n \geq n_0(M)$, \ as desired.
Hence for all \ $\varrho > \alpha$ \  and \ $y \geq \max\{\ty_0(n), \gamma n\}$ \ with
 \ $n \geq N_1$, \ we have
 \begin{align*}
  I_2 &=  \int_{\frac{qy}{a(n,y)}}^{qy} (\ee^{ht} - 1) \, F_{\eta_1}(\dd t)
       \leq \ee^{hqy} \PP\left(\eta_1 > \frac{qy}{a(n,y)}\right) \\
      &\leq \exp\left\{\frac{a(n,y)-K\varrho\log(a(n, y))}{K}\right\}
            B \left(\frac{q}{a(n,y)}\right)^{-\varrho} \PP(\eta_1 > y) \\
      &= B q^{-\varrho} \ee^{\frac{a(n,y)}{K}} \PP(\eta_1 > y)
       = B q^{-\varrho} (n \PP(\eta_1 > y))^{-\frac{1}{K}} \PP(\eta_1 > y) ,
 \end{align*}
 where we used that \ $1 < a(n, y) = -\log(n \PP(\eta_1 > y))$.
\ Putting together the bounds for \ $I_1$ \ and \ $I_2$ \ and using that
 \ $\frac{hqy}{a(n,y)} \leq \frac{1}{K}$ \ for \ $y \geq \gamma n$ \ with \ $n \geq N_1$, \ we
 obtain that
 \begin{equation}\label{expeq}
  \frac{\PP(\tS_n>y)}{n\PP(\eta_1>y)}
  \leq \exp\left\{n h \EE(\eta_1) \ee^{\frac{1}{K}}
                  + B q^{-\varrho} (n \PP(\eta_1 > y))^{1-\frac{1}{K}} - h y + a(n, y)\right\}
 \end{equation}
 for \ $y \geq \max\{\ty_0(n), \gamma n\}$ \ with \ $n \geq N_1$.
\ Noting that \ $n\PP(\eta_1 > y) \to 0$ \ uniformly for \ $y \geq \gamma n$ \ as \ $n \to \infty$
 \ (see, \eqref{help20}), we obtain that there exists \ $C_2 \in \RR_{++}$ \ such that the
 right-hand side of \eqref{expeq} can be bounded by
 \begin{align*}
  &C_2 \exp\left\{n h \EE(\eta_1) \ee^{\frac{1}{K}} - h y + a(n, y)\right\} \\
  &= C_2
     \exp\left\{h y \Biggl(\frac{\ee^{\frac{1}{K}}n\EE(\eta_1)}{y} - 1\Biggr) + a(n, y)\right\} \\
  &\leq C_2
        \exp\left\{\frac{a(n,y)-K\varrho\log(a(n,y))}{Kq}
                   \Biggl(\frac{\ee^{\frac{1}{K}}\EE(\eta_1)}{\gamma} - 1\Biggr) + a(n, y)\right\}
 \end{align*}
 for all \ $\varrho > \alpha$, \ sufficiently large \ $n \in \NN$ \ (greater than \ $N_1$) \ and
 \ $y \geq \max\{\ty_0(n), \gamma n\}$.
\ Since $\gamma > \EE(\eta_1)$, \ we can choose \ $K > 1$ \ sufficiently large such that
 \ $\frac{1}{\gamma} \ee^{\frac{1}{K}}\EE(\eta_1) < 1$, \ then we choose \ $q > 0$ \ sufficiently
 small such that
 \[
   \frac{1}{Kq} \biggl(\frac{\ee^{\frac{1}{K}}\EE(\eta_1)}{\gamma} - 1\biggr) + 2 < 0 ,
 \]
 i.e., \ $q < \frac{1}{2K} (1 - \frac{1}{\gamma} \ee^{\frac{1}{K}} \EE(\eta_1))$.
\ Then we have
 \begin{align*}
  \frac{\PP(\widetilde S_n>y)}{n\PP(\eta_1>y)}
  &\leq  C_2 \exp\left\{(a(n, y) - K \varrho \log(a(n, y))) (-2) + a(n, y)\right\} \\
  &= C_2 \exp\left\{2 K \varrho \log(a(n,y)) - a(n, y)\right\}
 \end{align*}
 for all \ $\varrho > \alpha$, \ sufficiently large \ $n \in \NN$ \ (greater than \ $N_1$) \ and
 \ $y \geq \max\{\ty_0(n), \gamma n\}$, \ where we used that
 \ $a(n, y) - K \varrho \log(a(n, y)) > 0$ \ for \ $y \geq \gamma n$ \ with \ $n \geq N_1$.
\ Here \ $C_2 \exp\left\{2 K \varrho \log(a(n, y)) - a(n, y)\right\}$ \ tends to \ $0$ \ uniformly
 for \ $y \geq \gamma n$ \ as \ $n \to \infty$, \ i.e.,
 \[
   \sup_{y\geq\gamma n} \exp\left\{2 K \varrho \log(a(n, y)) - a(n, y)\right\}
   = \exp\left\{\sup_{y\geq\gamma n} (2 K \varrho \log(a(n,y)) - a(n, y))\right\} \to 0
 \]
 as \ $n \to \infty$.
\ Indeed, this will be a consequence of
 \ $\sup_{y\geq\gamma n} (2 K \varrho \log(a(n, y)) - a(n, y)) \to - \infty$ \ as \ $n \to \infty$.
\ We have
 \begin{equation}\label{supsup}
  \sup_{y\geq\gamma n} (2 K \varrho \log(a(n, y)) - a(n, y)) \leq S_1(n) + S_2(n) ,
 \end{equation}
 where
 \begin{align*}
  S_1(n) &:= \sup_{y\geq\gamma n} \biggl(2 K \varrho \log(a(n,y)) - \frac{1}{2} a(n,y)\biggr) , \\
  S_2(n) &:= \sup_{y\geq\gamma n} \biggl(- \frac{1}{2} a(n,y)\biggr)
          = - \frac{1}{2} \inf_{y\geq\gamma n} a(n,y)
          \to - \infty \qquad \text{as \ $n \to \infty$.}
 \end{align*}
Moreover, \ $\lim_{x\to\infty} \frac{\log(x)}{x} = 0$ \ implies that there exists \ $\tM > 0$
 \ such that \ $\frac{\log(x)}{x} < \frac{1}{4K\varrho}$ \ for all \ $x \geq \tM$.
\ Since \ $\lim_{n\to\infty} \inf_{y\geq\gamma n} a(n,y) = \infty$, \ there exists
 \ $n_0(\tM) \in \NN$ \ such that \ $a(n, y) \geq \tM$ \ for all \ $y \geq \gamma n$ \ with
 \ $n \geq n_0(\tM)$.
\ Hence \ $\frac{\log(a(n,y))}{a(n,y)} < \frac{1}{4K\varrho}$ \ for all \ $y \geq \gamma n$ \ with
 \ $n \geq n_0(\tM)$, \ thus \ $2 K \varrho \log(a(n,y)) < \frac{1}{2} a(n,y)$.
\ Consequently, we obtain \ $S_1(n) \leq 0$ \ for all \ $n \geq n_0(\tM)$, \ and hence, by
 \eqref{supsup}, we conclude
 \ $\sup_{y\geq\gamma n} (2 K \varrho \log(a(n,y)) - a(n,y)) \to - \infty$ \ as \ $n \to \infty$,
 \ as desired.
So we have
 \[
   \lim_{n\to\infty} \sup_{y\geq \gamma n} \frac{\PP(\tS_n > y)}{n\PP(\eta_1 > y)} = 0.
 \]
Consequently, there exists an \ $N \in \NN$ \ such that
 \[
   \sup_{n\geq N,\,y\geq\gamma n} \frac{\PP(\tS_n>y)}{n\PP(\eta_1>y)} < \infty .
 \]
This, together with \eqref{help9_large} and \eqref{help10_large} yield that
 \begin{align}\label{help11_large}
  \sup_{n\geq N,\,y\geq\gamma n} \frac{\PP(S_n>y)}{n\PP(\eta_1>y)} < \infty .
 \end{align}
Finally, using the convolution property (see, Lemma \ref{Lem_conv}),
 \begin{align}\label{help12_large}
  \sup_{1\leq n\leq N,\,y\geq \gamma n} \frac{\PP(S_n>y)}{n\PP(\eta_1>y)}
  \leq \sum_{n=1}^N \sup_{y\geq \gamma n} \frac{\PP(S_n>y)}{n\PP(\eta_1>y)}
  < \infty .
 \end{align}
The desired statement readily follows from \eqref{help11_large} and \eqref{help12_large}.
\proofend

\section{Tail behavior of second-order Galton--Watson processes (without immigration) having
          regularly varying initial distributions}
\label{section_2GW}

\begin{Pro}\label{2GW_X_0_X_-1}
Let \ $(X_n)_{n\geq-1}$ \ be a second-order Galton--Watson process (without immigration) such that
 \ $X_0$ \ and \ $X_{-1}$ \ are independent, \ $X_0$ \ is regularly varying with index \ $\beta_0 \in \RR_+$, \ $X_{-1}$ \ is regularly
 varying with index \ $\beta_{-1} \in \RR_+$ \ and \ $m_\xi, m_\eta \in \RR_{++}$.
\ In case of \ $\max\{\beta_0, \beta_{-1}\} \in [1, \infty)$, \ assume additionally that there
 exists \ $r \in (\max\{\beta_0, \beta_{-1}\}, \infty)$ \ with \ $\EE(\xi^r) < \infty$ \ and
 \ $\EE(\eta^r) < \infty$.
\ Then for each \ $n \in \NN$,
 \[
   \PP(X_n > x)
   \sim \begin{cases}
         m_n^{\beta_0} \PP(X_0 > x) & \text{if \ $0 \leq \beta_0 < \beta_{-1}$,} \\
         m_n^{\beta_0} \PP(X_0 > x) + m_{n-1}^{\beta_{-1}} m_\eta^{\beta_{-1}} \PP(X_{-1} > x)
          & \text{if \ $\beta_0 = \beta_{-1}$,} \\
         m_{n-1}^{\beta_{-1}} m_\eta^{\beta_{-1}} \PP(X_{-1} > x)
          & \text{if \ $\beta_{-1} < \beta_0$}
        \end{cases}
 \]
 as \ $x \to \infty$, \ where \ $m_i$, \ $i \in \ZZ_+$, \ are given in Theorem \ref{Thm3} and hence,
 \ $X_n$ \ is regularly varying with index \ $\min\{\beta_0, \beta_{-1}\}$ \ for each \ $n \in \NN$.
\end{Pro}

\noindent{\bf First proof of Proposition \ref{2GW_X_0_X_-1}.}
Let us fix \ $n \in \NN$.
\ In view of the additive property \eqref{2GW_additive}, the independence of \ $X_0$ \ and \ $X_{-1}$, \
 and the convolution property of regularly varying distributions described in Lemma \ref{Lem_conv}, it is sufficient to prove
 \begin{equation}\label{ac}
  \PP\Biggl(\sum_{i=1}^{X_0} \zeta_{i,0}^{(n)} > x\Biggr) \sim m_n^{\beta_0} \PP(X_0 > x) , \qquad
  \PP\Biggl(\sum_{j=1}^{X_{-1}} \zeta_{j,-1}^{(n)} > x\Biggr)
  \sim m_{n-1}^{\beta_{-1}} m_\eta^{\beta_{-1}} \PP(X_{-1} > x)
 \end{equation}
 as \ $x \to \infty$.
\ These relations follow from Proposition \ref{FGAMSRS}, since
 \ $\EE(\zeta_{1,0}^{(n)}) =  m_n \in \RR_{++}$ \ and
 \ $\EE(\zeta_{1,-1}^{(n)}) = m_{n-1} m_\eta \in \RR_{++}$,
 \ $n \in \NN$, \ by \eqref{EXn}.
\proofend

\noindent{\bf Second proof of Proposition \ref{2GW_X_0_X_-1}.}
Let us fix \ $n \in \NN$.
\ In view of the additive property \eqref{2GW_additive}, the independence of \ $X_0$ \ and \ $X_{-1}$,
 \ and the convolution property of regularly varying distributions described in Lemma \ref{Lem_conv},
 it is sufficient to prove \eqref{ac}.
We show only the first relation in \eqref{ac}, since the second one can be proven in the same way.
Note that \ $\EE(\zeta_{1,0}^{(n)}) = m_n$ \ by \eqref{EXn}.
First, we prove
 \begin{align}\label{help3_liminf}
  \liminf_{x\to\infty}
   \frac{\PP\bigl(\sum_{i=1}^{X_0} \zeta_{i,0}^{(n)} > x\bigr)}{\PP(X_0 > x)}
  \geq m_n^{\beta_0} .
 \end{align}
Let \ $q \in (0, 1)$ \ be arbitrary.
For sufficiently large \ $x \in \RR_{++}$, \ we have
 \ $\lfloor(1+q)x/m_n\rfloor \geq 1$, \ since \ $m_n > 0$.
\ Using that for each \ $i \in \NN$, \ $\zeta_{i,0}^{(n)}$ \ is non-negative, we obtain
 \begin{align*}
  &\PP\Biggl(\sum_{i=1}^{X_0} \zeta_{i,0}^{(n)} > x\Biggr)
   \geq \sum_{k=\lfloor(1+q)x/m_n\rfloor}^\infty
         \PP\left(\sum_{i=1}^k \zeta_{i,0}^{(n)} > x\right) \PP(X_0 =k) \\
  &\geq \PP\Biggl(\sum_{i=1}^{\lfloor(1+q)x/m_n\rfloor} \zeta_{i,0}^{(n)} > x\Biggr)
        \sum_{k=\lfloor(1+q)x/m_n\rfloor}^\infty \PP(X_0 = k)\\
  &= \PP\Biggl(\frac{1}{\lfloor(1+q)x/m_n\rfloor}
               \sum_{i=1}^{\lfloor(1+q)x/m_n\rfloor} \zeta_{i,0}^{(n)}
               > \frac{x}{\lfloor(1+q)x/m_n\rfloor}\Biggr)
     \PP(X_0 \geq \lfloor(1+q)x/m_n\rfloor) \\
  &\geq \PP\Biggl(\frac{1}{\lfloor(1+q)x/m_n\rfloor}
                  \sum_{i=1}^{\lfloor(1+q)x/m_n\rfloor} \zeta_{i,0}^{(n)}
                  > \frac{x}{\lfloor(1+q)x/m_n\rfloor}\Biggr)
        \PP(X_0 > (1+q)x/m_n)
 \end{align*}
 for sufficiently large \ $x \in \RR_{++}$.
\ For sufficiently large \ $x \in \RR_{++}$, \ we have
 \ $\frac{x}{\lfloor(1+q)x/m_n\rfloor} \leq \frac{m_n}{1+(q/2)}$, \ since
 \ $\frac{x}{\lfloor(1+q)x/m_n\rfloor} \to \frac{m_n}{1+q}$ \ as \ $x \to \infty$ \ and
 \ $\frac{m_n}{1+q} < \frac{m_n}{1+(q/2)}$.
\ Hence, for sufficiently large \ $x \in \RR_{++}$, \ we have
 \[
   \PP\Biggl(\sum_{i=1}^{X_0} \zeta_{i,0}^{(n)} > x\Biggr)
   \geq \PP\Biggl(\frac{1}{\lfloor(1+q)x/m_n\rfloor}
                  \sum_{i=1}^{\lfloor(1+q)x/m_n\rfloor} \zeta_{i,0}^{(n)}
                  > \frac{m_n}{1+(q/2)}\Biggr)
        \PP\left(X_0 > \frac{(1+q)x}{m_n}\right) .
 \]
We have
 \begin{equation}\label{WLLN}
  \frac{1}{N} \sum_{i=1}^N \zeta_{i,0}^{(n)} {\as \EE(\zeta_{1,0}^{(n)})=}m_n \qquad
  \text{as \ $N \to \infty$}
 \end{equation}
 by the strong law of large numbers, hence \ $\frac{m_n}{1+(q/2)} < m_n$ \ yields
 \[
   \PP\Biggl(\frac{1}{\lfloor(1+q)x/m_n\rfloor}
             \sum_{i=1}^{\lfloor(1+q)x/m_n\rfloor} \zeta_{i,0}^{(n)}
             > \frac{m_n}{1+(q/2)}\Biggr)
   \to 1 \qquad \text{as \ $x \to \infty$.}
 \]
Thus, using that \ $X_0$ \ is regularly varying with index \ $\beta_0$, \ we have
 \begin{align*}
  \PP\Biggl(\frac{1}{\lfloor(1+q)x/m_n\rfloor}
            \sum_{i=1}^{\lfloor(1+q)x/m_n\rfloor} \zeta_{i,0}^{(n)}
            > \frac{m_n}{1+(q/2)}\Biggr)
  \PP\biggl(X_0 > \frac{(1+q)x}{m_n}\biggr)
  &\sim \PP\left(X_0 > \frac{(1+q)x}{m_n}\right) \\
  &\sim \Bigl(\frac{m_n}{1+q}\Bigr)^{\beta_0} \PP(X_0 > x)
 \end{align*}
 as \ $x \to \infty$.
\ Consequently,
 \[
   \liminf_{x\to\infty}
    \frac{\PP\bigl(\sum_{i=1}^{X_0} \zeta_{i,0}^{(n)} > x\bigr)}{\PP(X_0 > x)}
   \geq \Bigl(\frac{m_n}{1+q}\Bigr)^{\beta_0} , \qquad q \in (0, 1) ,
 \]
 and, by \ $q \downarrow 0$, \ we conclude \eqref{help3_liminf}.

Next, we prove
 \begin{align}\label{help3_limsup}
  \limsup_{x\to\infty} \frac{\PP\bigl(\sum_{i=1}^{X_0}\zeta_{i,0}^{(n)}>x\bigr)}{\PP(X_0>x)}
  \leq m_n^{\beta_0} .
 \end{align}
Let \ $q \in (0, 1)$ \ be arbitrary.
For sufficiently large \ $x \in \RR_{++}$, \  we have \ $\lfloor(1-q)x/m_n\rfloor\geq 1$, \ and
 hence
 \begin{align*}
  \PP\Biggl(\sum_{i=1}^{X_0} \zeta_{i,0}^{(n)} > x\Biggr)
  &\leq \PP(X_0 > \lfloor(1-q)x/m_n\rfloor)
        + \sum_{k=1}^{\lfloor(1-q)x/m_n\rfloor}
           \PP\left(\sum_{i=1}^k \zeta_{i,0}^{(n)} > x\right) \PP(X_0 =k) \\
  &= \PP\biggl(X_0 > \frac{(1-q)x}{m_n}\biggr)
     + \sum_{k=1}^{\lfloor(1-q)x/m_n\rfloor}
        \PP\left(\sum_{i=1}^k \zeta_{i,0}^{(n)} > x\right) \PP(X_0 =k) .
 \end{align*}
Since \ $X_0$ \ is regularly varying with index \ $\beta_0$, \ we have
 \[
  \PP\Bigl(X_0 > \frac{(1-q)x}{m_n}\Bigr)
  \sim \Bigl(\frac{m_n}{1-q}\Bigr)^{\beta_0} \PP(X_0 > x)
  \qquad \text{as \ $x \to \infty$,}
 \]
 hence, by taking the limit \ $q \downarrow 0$, \ we get \eqref{help3_limsup} provided we check
 \begin{equation}\label{p(x)}
  p(x, q)
  := \sum_{k=1}^{\lfloor(1-q)x/m_n\rfloor}
      \PP\left(\sum_{i=1}^k \zeta_{i,0}^{(n)} > x\right) \PP(X_0 =k)
  = \oo\left(\PP\left(X_0 > x\right)\right)
  \qquad \text{as \ $x\to\infty$}
 \end{equation}
 for all sufficiently small \ $q \in (0, 1)$.
\ (In fact, it will turn out that \eqref{p(x)} holds for any \ $q \in (0, 1)$.)

First, we consider the case \ $\beta_0 \in (0, 1)$.
\ Let \ $0 < \delta < (1-q)/m_n$.
\ Then for sufficiently large \ $x \in \RR_{++}$, \ we have
 \ $\lfloor\delta x\rfloor < \lfloor(1-q)x/m_n\rfloor$, \ and then
 \begin{align*}
  p(x, q) &= \sum_{k=1}^{\lfloor\delta x\rfloor}
              \PP\left(\sum_{i=1}^k \zeta_{i,0}^{(n)} > x\right) \PP(X_0 =k)
             + \sum_{k=\lfloor\delta x\rfloor+1}^{\lfloor(1-q)x/m_n\rfloor}
                \PP\left(\sum_{i=1}^k \zeta_{i,0}^{(n)} > x\right) \PP(X_0 =k) \\
          &=: p_1(x, \delta) + p_2(x, \delta, q) .
 \end{align*}
At first, we show that \ $p_2(x, \delta, q) = \oo(\PP(X_0 > x))$ \ as \ $x \to \infty$ \ for
 all \ $0 < \delta < (1-q)/m_n$.
\ Here, using that \ $\zeta_{i,0}^{(n)}$ \ is non-negative for each \ $i \in \NN$, \ we obtain
 \begin{align*}
  p_2(x, \delta, q)
  &\leq \PP\Biggl(\sum_{i=1}^{\lfloor(1-q)x/m_n\rfloor} \zeta_{i,0}^{(n)} > x\Biggr)
        \sum_{k=\lfloor\delta x\rfloor+1}^{\lfloor(1-q)x/m_n\rfloor} \PP(X_0 =k) \\
  &\leq \PP\Biggl(\sum_{i=1}^{\lfloor(1-q)x/m_n\rfloor} \zeta_{i,0}^{(n)} > x\Biggr)
        \PP(X_0 > \lfloor\delta x\rfloor)
   = \PP\Biggl(\sum_{i=1}^{\lfloor(1-q)x/m_n\rfloor} \zeta_{i,0}^{(n)} > x\Biggr)
     \PP(X_0 > \delta x) .
 \end{align*}
For sufficiently large \ $x \in \RR_{++}$, \ we have
 \ $\frac{x}{\lfloor(1-q)x/m_n\rfloor} \geq \frac{m_n}{1-(q/2)}$, \ since
 \ $\frac{x}{\lfloor(1-q)x/m_n\rfloor} \to \frac{m_n}{1-q}$ \ as \ $x \to \infty$ \ and
 \ $\frac{m_n}{1-q} > \frac{m_n}{1-(q/2)}$.
\ Hence, for sufficiently large \ $x \in \RR_{++}$, \ we have
 \begin{align*}
  \PP\Biggl(\sum_{i=1}^{\lfloor(1-q)x/m_n\rfloor} \zeta_{i,0}^{(n)} > x\Biggr)
  &= \PP\Biggl(\frac{1}{\lfloor(1-q)x/m_n\rfloor}
               \sum_{i=1}^{\lfloor(1-q)x/m_n\rfloor} \zeta_{i,0}^{(n)}
               > \frac{x}{\lfloor(1-q)x/m_n\rfloor}\Biggr) \\
  &\leq \PP\Biggl(\frac{1}{\lfloor(1-q)x/m_n\rfloor}
                  \sum_{i=1}^{\lfloor(1-q)x/m_n\rfloor} \zeta_{i,0}^{(n)}
                  > \frac{m_n}{1-(q/2)}\Biggr) .
 \end{align*}
Again by the strong law of large numbers (see \eqref{WLLN}), \ $\frac{m_n}{1-(q/2)} > m_n$ \ yields
 \[
   \PP\Biggl(\frac{1}{\lfloor(1-q)x/m_n\rfloor}
             \sum_{i=1}^{\lfloor(1-q)x/m_n\rfloor} \zeta_{i,0}^{(n)}
             > \frac{m_n}{1-(q/2)}\Biggr)
   \to 0 \qquad \text{as \ $x \to \infty$,}
 \]
 hence we obtain
 \begin{equation}\label{WLLN-}
  \PP\Biggl(\sum_{i=1}^{\lfloor(1-q)x/m_n\rfloor} \zeta_{i,0}^{(n)} > x\Biggr)
  \to 0 \qquad \text{as \ $x \to \infty$.}
 \end{equation}
Using that \ $X_0$ \ is regularly varying with index \ $\beta_0$, \ we have
 \ $\PP(X_0 > \delta x) \sim \delta^{-\beta_0} \PP(X_0 > x)$ \ as \ $x \to \infty$, \ hence
 \ $p_2(x, \delta, q) = \oo(\PP(X_0 > x))$ \ as \ $x \to \infty$ \ for all
 \ $0 < \delta < (1-q)/m_n$ \ and \ $q \in (0, 1)$.
\ Now we turn to prove
 \[
   \limsup_{\delta\downarrow 0} \limsup_{x\to\infty} \frac{p_1(x,\delta)}{\PP(X_0 > x)}   = 0 .
 \]
By Markov's inequality,
 \begin{align*}
  \PP\left(\sum_{i=1}^k \zeta_{i,0}^{(n)} > x\right)
  \leq \frac{1}{x} \sum_{i=1}^k \EE(\zeta_{i,0}^{(n)})
  = \frac{m_nk}{x}
 \end{align*}
 for all \ $k \in \NN$ \ and \ $x \in \RR_{++}$, \ and hence
 \begin{align*}
  p_1(x, \delta)
  &\leq \frac{m_n}{x}
        \sum_{k=0}^{\lfloor\delta x\rfloor} k\PP(X_0 =k)
   = \frac{m_n}{x} \EE(X_0 \bbone_{\{X_0\leq\lfloor\delta x\rfloor\}})
   = \frac{m_n}{x}
     \int_0^\infty \PP(X_0 \bbone_{\{X_0\leq\lfloor\delta x\rfloor\}} > t) \, \dd t \\
  &= \frac{m_n}{x}
     \int_0^{\lfloor\delta x\rfloor}
      \PP(X_0 \bbone_{\{X_0\leq\lfloor\delta x\rfloor\}} > t) \, \dd t
   \leq \frac{m_n}{x} \int_0^{\lfloor\delta x\rfloor} \PP(X_0 > t) \, \dd t .
 \end{align*}
Since \ $\RR_+ \ni x \mapsto \PP(X_0 > x)$ \ is locally integrable (due to the fact that it is
 bounded), it is integrable on intervals including \ $0$ \ as well, and since it is regularly
 varying (at infinity) with index \ $-\beta_0$, \ by Karamata's theorem (see Theorem
 \ref{Krthm}),
 \[
   \lim_{x\to\infty} \frac{x\PP(X_0>x)}{\int_0^x\PP(X_0>t)\,\dd t} = 1 - \beta_0,
 \]
 and hence
 \[
   \int_0^{\lfloor\delta x\rfloor} \PP(X_0 > t) \, \dd t
   \sim \frac{1}{1-\beta_0} \lfloor\delta x\rfloor \PP(X_0 > \lfloor\delta x\rfloor)
   = \frac{1}{1-\beta_0} \lfloor\delta x\rfloor \PP(X_0 > \delta x)
 \]
 as \ $x \to \infty$.
\ Then using that \ $\PP(X_0 > \delta x) \sim \delta^{-\beta_0} \PP(X_0 > x)$ \ as
 \ $x \to \infty$, \ we have
 \[
   \frac{p_1(x,\delta)}{\PP(X_0>x)}
   \leq \frac{m_n}{x} \frac{\int_0^{\lfloor\delta x\rfloor}\PP(X_0>t)\,\dd t}{\PP(X_0>x)}
   \sim \frac{m_n}{1-\beta_0} \delta \frac{\PP(X_0>\delta x)}{\PP(X_0>x)}
   \sim \frac{m_n}{1-\beta_0} \delta^{1-\beta_0}
 \]
 as \ $x \to \infty$.
\ Consequently,
 \[
   \limsup_{x\to\infty} \frac{p_1(x,\delta)}{\PP(X_0>x)}
   \leq \frac{m_n}{1-\beta_0} \delta^{1-\beta_0}
   \qquad \text{for all} \quad 0< \delta < \frac{1-q}{m_n} ,
 \]
 and hence
 \ $\limsup_{\delta\downarrow0} \limsup_{x\to\infty} \frac{p_1(x,\delta)}{\PP(X_0>x)}
    \leq \lim_{\delta\downarrow 0} \frac{m_n}{1-\beta_0} \delta^{1-\beta_0} = 0$.
\ Combining the parts we get \ $p(x, q) = \oo(\PP(X_0 > x))$ \ as \ $x \to \infty$ \ for any
 \ $q \in (0, 1)$, \ as desired.

Next, we consider the case \ $\beta_0 \in (1, 2)$.
\ Using Lemma \ref{Lem_L_const}, we check that there exists a non-negative random variable
 \ $\tzeta^{(n)}$ \ having the following properties:
 \begin{itemize}
  \item
   $\tzeta^{(n)}$ \ is regularly varying with index \ $\beta_0$,
  \item
   $\PP(\zeta_{1,0}^{(n)} > x) \leq \PP(\tzeta^{(n)} > x)$, \ $x \in \RR_+$,
  \item
   $\PP(\tzeta^{(n)} > x) = \oo(\PP(X_0 > x))$ \ as \ $x \to \infty$,
  \item
   $\EE(\zeta_{1,0}^{(n)}) \leq \EE(\tzeta^{(n)}) < \infty$.
 \end{itemize}
By Lemma \ref{Lem_seged_momentr}, \ $\EE((\zeta_{1,0}^{(n)})^r) < \infty$, \ and hence, by Lemma
 \ref{exposv}, \ $\PP(\zeta_{1,0}^{(n)} > x) = \oo(\PP(X_0 > x))$ \ as \ $x \to \infty$.
\ Thus, by Lemma \ref{Lem_L_const}, there exists a monotone increasing, right-continuous, slowly
 varying (at infinity) function \ $L_{\tzeta^{(n)}}$ \ such that \ $L_{\tzeta^{(n)}}(x) \geq 1$,
 \ $x \in \RR_+$, \ $\lim_{x\to\infty} L_{\tzeta^{(n)}}(x) = \infty$ \ and
 \ $\lim_{x\to\infty} L_{\tzeta^{(n)}}(x) \frac{\PP(\zeta_{1,0}^{(n)}>x)}{\PP(X_0>x)} = 0$.
\ Hence, using also that \ $\PP(X_0 \geq x) \leq 1$, \ $x \in \RR_+$, \ there exists
 \ $x' \in \RR_+$ \ such that
 \ $L_{\tzeta^{(n)}}(x) \frac{\PP(\zeta_{1,0}^{(n)}>x)}{\PP(X_0>x)} \leq 1$ \ and
 \ $\frac{\PP(X_0>x)}{L_{\tzeta^{(n)}}(x)}  \leq 1$ \ hold for all \ $x \geq x'$.
\ Let \ $\tzeta^{(n)}$ \ be a random variable such that
 \[
   \PP(\tzeta^{(n)} > x) := \begin{cases}
                             1 &  \text{if \ $x \leq x'$,} \\
                             \frac{\PP(X_0>x)}{L_{\tzeta^{(n)}}(x)} & \text{if \ $x > x'$.}
                            \end{cases}
 \]
Such a non-negative random variable exists, since
 \ $\RR_{++} \ni x \mapsto \frac{\PP(X_0>x)}{L_{\tzeta^{(n)}}(x)}$ \ is monotone decreasing,
 converges to \ $0$ \ as \ $x \to \infty$ \ and right-continuous.
For all \ $q \in \RR_{++}$,
 \[
   \lim_{x\to\infty} \frac{\PP(\tzeta^{(n)}>qx)}{\PP(\tzeta^{(n)}>x)}
   = \lim_{x\to\infty}
      \frac{L_{\tzeta^{(n)}}(x)}{L_{\tzeta^{(n)}}(qx)} \frac{\PP(X_0>qx)}{\PP(X_0>x)}
   = 1 \cdot q^{-\beta_0}
   = q^{-\beta_0} ,
 \]
 yielding that \ $\tzeta^{(n)}$ \ is regularly varying with index \ $\beta_0$.
\ For \ $x \leq x'$, \ we have \ $\PP(\zeta_{1,0}^{(n)} > x) \leq 1 = \PP(\tzeta^{(n)} > x)$.
\ For \ $x > x'$, \ we have
 \[
   \PP(\zeta_{1,0}^{(n)} > x)
   = L_{\tzeta^{(n)}}(x) \frac{\PP(\zeta_{1,0}^{(n)}>x)}{\PP(X_0>x)} \PP(\tzeta^{(n)} > x)
   \leq \PP(\tzeta^{(n)} > x) .
 \]
Further,
 \[
   \lim_{x\to\infty} \frac{\PP(\tzeta^{(n)}>x)}{\PP(X_0>x)}
   = \lim_{x\to\infty} \frac{\PP(X_0>x)}{L_{\tzeta^{(n)}}(x)\PP(X_0>x)}
   = 0 ,
 \]
 since \ $\lim_{x\to\infty} L_{\tzeta^{(n)}}(x) = \infty$.
\ Since \ $\PP(\zeta_{1,0}^{(n)} > x) \leq \PP(\tzeta^{(n)} > x)$, \ $x \in \RR_+$, \ we have
 \[
   \EE(\zeta_{1,0}^{(n)})
   = \int_0^\infty \PP(\zeta_{1,0}^{(n)} > x) \, \dd x
   \leq \int_0^\infty \PP(\tzeta^{(n)} > x) \, \dd x
   = \EE(\tzeta^{(n)}) ,
 \]
 and since \ $\tzeta^{(n)}$ \ is regularly varying with index \ $\beta_0 \in (1, 2)$, \ we
 have \ $\EE(\tzeta^{(n)}) < \infty$.

Let \ $(\tzeta_j^{(n)})_{j\in\NN}$ \ be a sequence of independent identically distributed random
 variables with common distribution as that of \ $\tzeta^{(n)}$.
\ By some properties of first order stochastic dominance (see, e.g., Shaked and Shanthikumar
 \cite[part (b) of Theorem 1.A.3 and Theorem 1.A.4]{ShaSha}), we have
 \begin{equation}\label{fod}
  \PP\left(\sum_{i=1}^k \zeta_{i,0}^{(n)} > x\right)
  \leq \PP\left(\sum_{i=1}^k \tzeta_i^{(n)} > x\right)
 \end{equation}
 for all \ $x \in \RR_+$ \ and \ $k \in \NN$.
\ Put \ $\tm_n := \EE(\tzeta^{(n)})$.
\ Let us consider the decomposition
 \begin{align*}
  p(x, q) &= \sum_{k=1}^{\lfloor(1-q)x/\tm_n\rfloor}
              \PP\left(\sum_{i=1}^k \zeta_{i,0}^{(n)} > x\right) \PP(X_0 =k) \\
          &\quad
             + \sum_{k=\lfloor(1-q)x/\tm_n\rfloor+1}^{\lfloor(1-q)x/m_n\rfloor}
                \PP\left(\sum_{i=1}^k \zeta_{i,0}^{(n)} > x\right) \PP(X_0 =k)
           =: p_1(x, q) + p_2(x, q) , \qquad x \in \RR_+ .
 \end{align*}
Here \ $m_n \leq \tm_n$, \ and hence \ $\lfloor(1-q)x/\tm_n\rfloor \leq \lfloor(1-q)x/m_n\rfloor$,
 \ $x \in \RR_+$, \ $q \in (0, 1)$.
\ Applying Theorem \ref{Thm_largedev} with \ $\gamma := \frac{\tm_n}{1-q} > \tm_n$, \ we conclude
 the existence of a  constant \ $C(q, n) \in \RR_{++}$ \ (not depending on \ $k$ \ and \ $x$, \ but
 on \ $q$ \ and \ $n$) \ such that
 \begin{equation}\label{largedev}
  \PP\left(\sum_{i=1}^k \tzeta_i^{(n)} > x\right) \leq C(q, n) k \PP(\tzeta^{(n)} > x) \qquad
  \text{for all \ $x \geq \gamma k$, \  $k\in\NN$.}
 \end{equation}
Using \eqref{fod} and \eqref{largedev}, we obtain
 \begin{align*}
  p_1(x, q)
  &\leq \sum_{k=1}^{\lfloor(1-q)x/\tm_n\rfloor}
         \PP\left(\sum_{i=1}^k \tzeta_i^{(n)} > x\right) \PP(X_0 =k) \\
  &\leq C(q, n) \sum_{k=1}^{\lfloor(1-q)x/\tm_n\rfloor} k \PP(\tzeta^{(n)} > x) \PP(X_0 =k)
   \leq C(q, n) \EE(X_0) \PP(\tzeta^{(n)} > x) ,
   \qquad x \in \RR_+ .
 \end{align*}
Hence for each \ $q \in (0, 1)$,
 \[
   \limsup_{x\to\infty} \frac{p_1(x,q)}{\PP(X_0>x)}
   \leq C(q, n) \EE(X_0) \limsup_{x\to\infty} \frac{\PP(\tzeta^{(n)}>x)}{\PP(X_0>x)} = 0 ,
 \]
 where the last step follows by the corresponding property of \ $\tzeta^{(n)}$.
\ Moreover,
 \begin{align*}
  p_2(x, q)
  &\leq \PP\left(\sum_{i=1}^{\lfloor(1-q)x/m_n\rfloor} \zeta_{i,0}^{(n)} > x\right)
        \sum_{k=\lfloor(1-q)x/\tm_n\rfloor+1}^{\lfloor(1-q)x/m_n\rfloor} \PP(X_0 =k) \\
  &\leq \PP\left(\sum_{i=1}^{\lfloor(1-q)x/m_n\rfloor} \zeta_{i,0}^{(n)} > x\right)
        \PP(X_0 > \lfloor(1-q)x/\tm_n\rfloor) \\
  &= \PP\left(\sum_{i=1}^{\lfloor(1-q)x/m_n\rfloor} \zeta_{i,0}^{(n)} > x\right)
     \PP\left(X_0 > \frac{(1-q)x}{\tm_n}\right) .
 \end{align*}
Since \ $X_0$ \ is regularly varying with index \ $\beta_0$, \ we have
 \[
   \lim_{x\to\infty}
    \frac{\PP\left(X_0>\frac{(1-q)x}{\tm_n}\right)}{\PP(X_0>x)}
   = \left(\frac{\tm_n}{1-q}\right)^{\beta_0} ,
 \]
 hence, for each \ $q \in (0, 1)$, \ applying \eqref{WLLN-}, we conclude
 \[
   \limsup_{x\to\infty} \frac{p_2(x,q)}{\PP(X_0>x)}
   \leq 0 \cdot \left(\frac{\tm_n}{1-q}\right)^{\beta_0} = 0 .
 \]
 Finally, we turn to the case \ $\beta_0 = 1$.
\ For each \ $q \in (0, 1)$, \ we have
 \begin{align*}
  p(x, q) &= \sum_{k=1}^{\lfloor(1-q)x/m_n\rfloor}
              \PP\left(\sum_{i=1}^k \zeta_{i,0}^{(n)} > x\right) \PP(X_0 =k) \\
          &= \sum_{k=1}^{\lfloor(1-q)x/m_n\rfloor}
              \PP\left(\sum_{i=1}^k \zeta_{i,0}^{(n)}- k m_n > x - k m_n\right) \PP(X_0 = k) .
 \end{align*}
Let \ $r' \in (1, 2]$.
\ According to Lemma 2.1 in Robert and Segers \cite{RobSeg} with \ $\gamma = \frac{m_nq}{1-q}$,
 \ there exist positive numbers \ $v$ \ and \ $C = C(v, q, n)$ \ such that for all \ $x \in \RR_+$
 \ and \ $k \in \NN$ \ with \ $k \leq \lfloor(1-q)x/m_n\rfloor$,
 \[
   \PP\left(\sum_{i=1}^k \zeta_{i,0}^{(n)}- km_n > x-km_n\right)
   \leq k \PP\bigl(\zeta_{1,0}^{(n)} - m_n > v(x-km_n)\bigr) + \frac{C}{(x-km_n)^{r'}} ,
 \]
 since \ $x - k m_n \geq \gamma k$ \ for all \ $x \in \RR_+$ \ and \ $k \in \NN$ \ with
 \ $k \leq \lfloor(1-q)x/m_n\rfloor$.
\ Consequently,
 \begin{align*}
  p(x, q)
  &= \sum_{k=1}^{\lfloor(1-q)x/m_n\rfloor}
      \PP\left(\sum_{i=1}^k \zeta_{i,0}^{(n)}- k m_n > x - k m_n\right) \PP(X_0 = k) \\
  &\leq \sum_{k=1}^{\lfloor(1-q)x/m_n\rfloor}
         \PP(X_0 = k)
         \left(k \PP\bigl(\zeta_{1,0}^{(n)} - m_n > v(x - k m_n)\bigr)
               + \frac{C}{(x-km_n)^{r'}}\right) \\
  &\leq \sum_{k=1}^{\lfloor(1-q)x/m_n\rfloor}
         \PP(X_0 = k)
         \left(k \PP\bigl(\zeta_{1,0}^{(n)} > v(x - k m_n)\bigr)
               + \frac{C}{(x-km_n)^{r'}}\right) \\
  &\leq \PP(\zeta_{1,0}^{(n)} > q v x) \sum_{k=1}^{\lfloor(1-q)x/m_n\rfloor} k \PP(X_0 = k)
        + \frac{C}{(qx)^{r'}} \\
  &\leq \PP(\zeta_{1,0}^{(n)} > q v x)
        \EE\big(X_0 \bbone_{\{X_0\leq\lfloor(1-q)x/m_n\rfloor\}}\big)
        + \frac{C}{(qx)^{r'}} ,
 \end{align*}
 where for the last but one step, we used that \ $x - k m_n \geq q x$ \ for
 \ $k \in \{1, \ldots, \lfloor(1-q)x/m_n\rfloor$\}.
\ Since \ $r' \in (1, 2]$, \ by Lemma \ref{sv}, we have \ $C/(qx)^{r'} = \oo(\PP(X_0 >x))$ \ as
 \ $x \to \infty$, \ so we only have to work with the first term.
If \ $\EE(X_0) < \infty$, \ then
 \ $\EE\big(X_0 \bbone_{\{X_0\leq\lfloor(1-q)x/m_n\rfloor\}}\big) \leq \EE(X_0) < \infty$ \ also
 holds, and
 \[
   \frac{\PP(\zeta_{1,0}^{(n)}>qvx)}{\PP(X_0>x)}
   = \frac{\PP(X_0>qvx)}{\PP(X_0>x)} \cdot \frac{\PP(\zeta_{1,0}^{(n)}>qvx)}{\PP(X_0>qvx)}
   \to 0 \qquad \text{as $x \to \infty$,}
 \]
 where we used that \ $X_0$ \ is regularly varying with index \ $1$, \ and that
 \ $\PP(\zeta_{1,0}^{(n)} > x) = \oo(\PP(X_0 > x))$ \ as \ $x \to \infty$ \ also holds (as it was
 already proved earlier).
Now we consider the case \ $\EE(X_0) = \infty$.
\ By Markov's inequality,
 \ $\PP(\zeta_{1,0}^{(n)} > q v x) \leq \EE((\zeta_{1,0}^{(n)})^r)/(qvx)^r$
 \ (note that in this case, \ $\EE((\zeta_{1,0}^{(n)})^r)$ \ exists, see Lemma
 \ref{Lem_seged_momentr}), and using the fact that
 \ $\limsup_{x\to\infty} \frac{\EE(X_0\bbone_{\{X_0\leq x\}})}{x^s\PP(X_0 >x)} = 0$ \ for
 some \ $1 < s < r$ \ (see the remark after Theorem 3.2 in Robert and Segers \cite{RobSeg}), we
 have
 \begin{align*}
  &\frac{\PP(\zeta_{1,0}^{(n)}>qvx)
         \EE\big(X_0\bbone_{\{X_0\leq\lfloor(1-q)x/m_n\rfloor\}}\big)}
        {\PP(X_0>x)}
   \leq \frac{\EE((\zeta_{1,0}^{(n)})^r)
              \EE\big(X_0\bbone_{\{X_0\leq\lfloor(1-q)x/m_n\rfloor\}}\big)}
             {(qvx)^r\PP(X_0>x)} \\
  &= \frac{\EE((\zeta_{1,0}^{(n)})^r)}{(qv)^r}
     \cdot \frac{\EE\big(X_0\bbone_{\{X_0\leq\lfloor(1-q)x/m_n\rfloor\}}\big)}{x^s\PP(X_0>x)}
     \cdot \frac{1}{x^{r-s}} \\
  &= \frac{\EE((\zeta_{1,0}^{(n)})^r)}{(qv)^r} \cdot \frac{\lfloor(1-q)x/m_n\rfloor^s}{x^s}
     \cdot \frac{\PP\left(X_0>\lfloor(1-q)x/m_n\rfloor\right)}{\PP(X_0>x)} \\
  &\phantom{=\;}
     \times \frac{\EE\big(X_0\bbone_{\{X_0\leq\lfloor(1-q)x/m_n\rfloor\}}\big)}
                {\lfloor(1-q)x/m_n\rfloor^s\PP(X_0>\lfloor(1-q)x/m_n\rfloor)}
     \cdot \frac{1}{x^{r-s}}
   \to 0 \qquad \text{as \ $x \to \infty$.}
 \end{align*}
Putting parts together, we have \ $p(x, q) = \oo(\PP(X_0 > x))$ \ as \ $x \to \infty$, \ as
 desired.
\proofend

\begin{Rem}
For a corresponding result for (first-order) Galton--Watson processes (without immigration), see
 Barczy et al.\ \cite[Proposition 2.2]{BarBosPap}.
A formal application of Proposition \ref{2GW_X_0_X_-1} also gives this result, namely, for each
 \ $n \in \NN$, \ we have \ $\PP(X_n > x) \sim m_\xi^{n\beta_0} \PP(X_0 > x)$ \ as
 \ $x \to \infty$.
\ In case of \ $m_\xi = 0$ \ and \ $m_\eta \in \RR_{++}$, \ Proposition 2.2 in Barczy et al.\
 \cite{BarBosPap} gives that \ $\PP(X_n > x) \sim m_\eta^{\frac{n}{2}\beta_{0}} \PP(X_0 > x)$ \ as
 \ $x \to \infty$ \ if \ $n \in \NN$ \ is even, and
 \ $\PP(X_n > x) \sim m_\eta^{\frac{n+1}{2}\beta_{-1}} \PP(X_{-1} > x)$ \ as \ $x \to \infty$ \ if
 \ $n \in \NN$ \ is odd.
\proofend
\end{Rem}

\section*{Acknowledgements}
We would like to thank the referee and Prof.\ Yuliya Mishura, Co-editor-in-chief, for their comments that helped us to improve the paper.

\end{document}